\documentclass[useAMS,usenatbib]{mnras}
\usepackage{amsmath}
\usepackage{amssymb}
\usepackage{graphicx,txfonts,fleqn,enumerate}
\usepackage{color}        
\newcommand   {\change}[1] {#1}   

\renewcommand*{\emph}[1] {\textit{#1}}
\renewcommand*{\vec}[1]  {\boldsymbol{#1}}
\newcommand*  {\diff}    {\mathop{}\!\mathrm{d}}
\newcommand*  {\Exp}[1]  {\mathrm{e}^{#1}}
\newcommand*  {\op}[1]   {{\hat{#1}}}
\newcommand*  {\map}[2]  {\Exp{{#1}\op{#2}}}
\newcommand*  {\altmap}[2] {\varphi_{#1}^{#2}}
\newcommand*  {\bV}      {\bar{V}}
\newcommand*  {\bVx}     {\bar{V}_{n-1}}
\renewcommand*{\ij}      {{i\!j}}
\newcommand*  {\pdiff}[2]{\frac{\partial{#1}}{\partial{#2}}}
\newcommand*  {\sub}[2]  {{#1}_{\mathrm{#2}}}
\newcommand*  {\opn}		{\change{\{}}
\newcommand*  {\cls}		{\change{\}}}

\title[Symplectic fourth-order $N$-body maps]
	  {\boldmath Symplectic fourth-order maps for the collisional $N$-body problem}
\date{Accepted . Received ; }
\author[Walter Dehnen \& David M.\ Hernandez]
	   {Walter Dehnen$^1$\thanks{wd11@le.ac.uk,\;dmhernan@mit.edu} and David M.\ Hernandez$^2$\footnotemark[1]
	   \\
	   $^1$Department of Physics \& Astronomy, University of Leicester,
	   	Leicester, LE1 7RH, UK\\
	   $^2$Department of Physics and Kavli Institute for Astrophysics and Space Research, Massachusetts Institute of Technology, 77 Massachusetts Ave.,\\
	   \phantom{$^2$}Cambridge, Massachusetts 02139, USA
	   \\
	   \vspace*{-2.5ex}
	   }
	
\pagerange{\pageref{firstpage}--\pageref{lastpage}}
\pubyear{2016}
\begin{document}
\label{firstpage}
\maketitle
\begin{abstract}
We study analytically and experimentally \change{certain} symplectic and time-reversible $N$-body integrators which employ a Kepler solver for each pair-wise interaction, \change{including} the method of \citeauthor{HernandezBertschinger2015}. Owing to the Kepler solver, these methods treat close two-body interactions correctly, while close three-body encounters contribute to the truncation error at second order and above. The second-order errors can be corrected to obtain a fourth-order scheme with little computational overhead. We generalise this map to an integrator which employs a Kepler solver only for selected interactions and yet retains fourth-order accuracy without backward steps. In this case, however, two-body encounters not treated via a Kepler solver contribute to the truncation error.
\end{abstract}

\begin{keywords}
gravitation - methods: analytical - methods: numerical - celestial mechanics - globular clusters: general - planets and satellites: dynamical evolution and stability\vspace*{-1ex}
\end{keywords}

\section{Introduction}
The gravitational $N$-body problem has been studied ever since Newton first wrote down his universal gravitational law of attraction. The $N$-body problem appears often in dynamical astronomy, for example  planetary systems, stellar associations, star clusters, galaxies, dark matter haloes, and even the universe as a whole can be modelled to good approximation as $N$-body problems \citep{HeggieHut2003}, although other, typically less accurate, alternative models are possible in some cases. No analytic solutions to the $N$-body problem exist for $N>2$, except for few cases \change{without practical relevance,} such as the five families of solutions found by \cite{Euler1767} and \cite{Lagrange1772}, and numerical integration is required instead.

If the $N$-body method is used to model a collision-less system (where two-body encounters are dynamically unimportant), encounters between the simulation particles introduce relaxation into the model not present in the actual system. These artificial effects can be reduced \change{(but not eliminated)} by softening the gravitational inter-particle forces at small distances \change{\citep{DehnenRead2011}}, which in turn significantly simplifies the $N$-body dynamics and allows the use of comparatively simple integration techniques, such as the leapfrog integrator \citep{Stoermer1907, Verlet1967}\footnote{The leapfrog integrator has been independently discovered several times, and \change{was implicitly used by \cite[][figure for theorem I in book I]{Newton1687} as later discovered by Verlet himself \citep*{HairerLubichWanner2006}}.}.

Here, we are instead concerned with the collisional $N$-body problem, which emerges for example when modelling the planetary systems including our own, planetesimals in a circum-stellar disc, or a globular cluster. In this case, the accurate long-term time integration of the unsoftened gravitational forces poses a formidable problem. Here `long term' means several Lyapunov times or when a conventional integrator becomes unreliable due to accumulation of truncation errors, whichever is shorter. A major problem arises from the dynamical stiffness of these systems in the sense that the relevant time scales differ by orders of magnitude: already a simple elliptic or hyperbolic orbit poses problems for numerical integration owing to the large variation of angular speed, i.e.\ of the local orbital time scale.

Since the $N$-body problem comprises a Hamiltonian system, symplectic, or more broadly geometric, numerical time integration\footnote{A symplectic integrator advances the system by a canonical map which is close to that of the actual Hamiltonian. As a consequence, the geometric structure of phase space and the Poincar\'{e} invariants are exactly preserved. Many symplectic integrators also \change{exactly} conserve all first integrals, except for the Hamiltonian, which tends to have bounded error.} \citep*{HairerLubichWanner2006} provides a useful framework for the $N$-body problem. Unfortunately, symplectic integration has not been widely implemented for the study of the collisional $N$-body problem. Switching methods, which change between different symplectic integrators, have been proposed for the of study single-star planetary systems \citep*{Chambers1999, KvaernoLeimkuhler2000, DuncanLevisonLee1998}. Unfortunately, tests indicate these methods may break time-reversibility and symplecticity \citep{Hernandez2016}. Another possibility to deal with the varying time-scales is to transform to another time variable \citep[Sundman transform, see][]{LeimkuhlerReich2004} and apply a symplectic method in the resulting extended system \citep{PretoTremaine1999, MikkolaTanikawa1999}, \change{but} such methods cannot be efficient for $N\gg2$.

An alternative to exactly symplectic integrators are time-reversible geometric integration methods, which share many desirable properties with symplectic integrators \citep{HairerLubichWanner2006}. When modelling globular clusters, a common such  integration method is the implicit fourth-order Hermite integrator \citep{Makino1991}, which requires an iterative solution, but in practice often only one iteration is used, violating exact time symmetry. Even when iterating to convergence, the efficient adaptation of individual discrete step sizes cannot be reconciled with exact time symmetry (Dehnen 2016, in preparation). \cite*{KokuboYoshinagaMakino1998} argue that this is tolerable if only few step-size changes occur, such as in planetary systems with only near-circular orbits. \cite*{HutMakinoMcMillan1995} proposed a symmetrisation procedure for any integrator and \cite{MakinoEtAl2006} extended this procedure to adaptation of individual particle step sizes. However, the resulting method involves the solution of a large implicit system of equations requiring an excessive amount of computational effort and has not been used in practice.

Because of these complications, contemporary methods for the integration of planetary systems employ a fixed global time step. A recent progress was the introduction of a symplectic and time-reversible map which treats close two-body encounters exactly \citep{HernandezBertschinger2015} and is efficient for planetary-system integration \citep{Hernandez2016}. In this present study, we show that the integrator \change{of \citeauthor{HernandezBertschinger2015}} is \change{still only} second-order accurate, but can be made fourth-order accurate with relatively little additional computational effort. We also discuss the option to treat only selected pair-wise interactions exactly (to improve efficiency) and yet keep the overall integration accuracy at fourth order.

This paper is organised as follows. Section~\ref{sec:prelim} reviews background concepts on symplectic integration and re-analyses the popular leapfrog (St{\o}rmer-Verlet) integrator, Section~\ref{sec:kepler} discusses \change{the integrator of \citeauthor{HernandezBertschinger2015}}, introduces its fourth-order extension, and presents some numerical tests. Integrators which use a Kepler solver selectively are considered in Section~\ref{sec:mix}, including our novel fourth-order hybrid integrator. The appendices provide some detailed  calculations and discuss implementation details. 

\section{\change{Symplectic maps from operator splitting}}
\label{sec:prelim}
The time-evolution for systems with Hamiltonian function $H$ is a continuous canonical transformation governed by
\begin{equation}\label{eq:time:evolve}
	\frac{\diff{w}}{\diff{t}} = \op{H}w \equiv \{w,H\}
\end{equation}
with $\{,\}$ the Poisson bracket \change{and $w\equiv\{\vec{x}_i,\vec{p}_i\}$ the set of all coordinates and momenta}. This equation defines the operator $\op{H}$, also known as \emph{Lie operator} of the function $H$ \citep{DragtFinn1976}, and has formal solution
\begin{equation} \label{eq:map:def}
	w(t+h) = \map{h}{H}  w(t)
\end{equation}
If no exact solution to~(\ref{eq:time:evolve}) exists, the time-evolution operator $\map{h}{H}$ has no finite expression, and instead a numerical solution is required. A symplectic integrator is such a numerical method that preserves the symplecticity (canonical nature) of the map $\map{h}{H}$. If one can split $H=A+B$ such that equation~(\ref{eq:time:evolve}) with $H$ replaced by $A$ or $B$ can be solved exactly, then a symplectic integrator can be constructed as composition of the maps $\map{h}{A}$ and $\map{h}{B}$. The simplest such method is the symplectic Euler method\begin{equation}\label{eq:Forward:Euler}
	\map{h}{H} \to \map{h}{A}\map{h}{B}.
\end{equation}
Thus, this method applies the time evolution due to $B$ followed by that due to $A$. The error made by the symplectic Euler method can be expressed by the \cite{Campbell1996, Campbell1997}-\cite{Baker1902b, Baker1905}-\cite{Hausdorff1906} formula \citep{Dynkin1947}
\begin{eqnarray}\label{eq:BCH}
	\log\big(\Exp{X}\,\Exp{Y}\big) &=& 
	X+Y + \tfrac{1}{2}[X,Y] + \tfrac{1}{12}\big([X,\change{[}X,Y\change{]}]
				 + [Y,\change{[}Y,X\change{]}]\big)
	\;\dots
\end{eqnarray}
with $[X,Y]\equiv X Y-Y X$ the usual commutator. Using the Jacobi identity
\begin{equation}
	\{\opn A,B\cls,C\} + \{\opn B,C\cls,A\} + \{\opn C,A\cls,B\} = 0,
\end{equation}
it can be shown that the Lie operator of the Poisson bracket $\{A,B\}$ of two phase-space functions \change{$A$ and $B$} is the commutator $[\op{B},\op{A}]$ of their Lie operators
\begin{equation}\label{eq:CommutePoisson}
\change{
	\widehat{\{B,A\}} = \{.,\{B,A\}\} = [\op{A},\op{B}]
}
\end{equation}
\change{which can be applied recursively: $\widehat{\{\{C,B\},A\}}=[\op{A},[\op{B},\op{C}]]$ etc.} Together with the distributive property \change{$\op{A}+\op{B}=\widehat{A+B}$} and the Camp\-bell-Baker-Haussdorff formula~(\ref{eq:BCH}) this implies that, under some conditions described below, the symplectic Euler method~(\ref{eq:Forward:Euler}) actually evolves the system under the surrogate Hamiltonian $\tilde{H} = H + \sub{H}{err}(h)$ with
\begin{equation}\label{eq:Forward:Euler:Hsurr}
	\sub{H}{err} =
	\frac{h}{2} \{B,A\} +
	\frac{h^2}{12} \{\opn B,A\cls,A\} +
	\frac{h^2}{12} \{\opn A,B\cls,B\} +
	\mathcal{O}(h^3),
\end{equation}
i.e.\ makes an error $\mathcal{O}(h^2)$ in the coordinates per time step and $\mathcal{O}(h)$ in the energy. A better integrator is the leapfrog or Verlet method
\begin{equation}\label{eq:Leapfrog}
	\map{h}{H} \to \Exp{\frac{h}{2}\op{A}}\,\map{h}{B}\,\Exp{\frac{h}{2}\op{A}}.
\end{equation}
Applying equation~(\ref{eq:Forward:Euler:Hsurr}) twice, we find for the leapfrog
\begin{equation} \label{eq:Leapfrog:Herr}
	\sub{H}{err}=
		-\frac{h^2}{24}\{\{B,A\},A\}
		+\frac{h^2}{12}\{\{A,B\},B\}
		+\mathcal{O}(h^4).
\end{equation}
In particular, \change{no odd powers of $h$ appear}, which is true for any self-adjoint integrator\footnote{\change{If $\altmap{h}{-1}$ is t}he inverse  of a phase-space map $\altmap{h}{}$\change{,} defined by the condition that the composite $\altmap{h}{}\altmap{h}{-1}$ is the identity map\change{, then} $\altmap{h}{\dag}\equiv\altmap{-h}{-1}$ is called \emph{adjoint} to $\altmap{h}{{}}$. For self-adjoint maps $\altmap{h}{\dag}=\altmap{h}{{}}$, which implies $\altmap{-h}{{}}=\altmap{h}{-1}$, i.e.\ these maps are also reversible or \change{time} symmetric.} for symmetry reasons.  

$\tilde{H}$ is a power series in $h$ that can converge or diverge. \change{In case of convergence}, $\tilde{H}$ is conserved and has properties of a Hamiltonian \citep{DragtFinn1976}. We have \change{never} found evidence \change{for divergence whenever we tested it, but} addressing \change{this} issue \change{further} is beyond the scope of this paper. \change{Instead,} we generally assume $\tilde{H}$ converges as is often done in the literature.
 
\subsection{The Leapfrog \boldmath $N$-body integrator}
The traditional splitting of the $N$-body Hamiltonian is in kinetic and potential energies,
\begin{equation} \label{eq:T,V}
	T = \sum_i \frac{\vec{p}_i^2}{2m_i},\quad
	V = \sum_{i,j<i} V_{\ij} = \frac{1}{2} \sum_{i,j} V_{\ij}
	\quad\text{with}\quad
	V_{\ij} = -\frac{Gm_im_{\!j}}{|\vec{x}_{\ij}|},
\end{equation}
where $\vec{x}_{\ij}\equiv\vec{x}_i-\vec{x}_{\!j}$ is the distance vector. The map $\map{h}{T}$ obtains a simple drift at constant velocity and $\map{h}{V}$ a kick, a change of velocity at fixed position. There are two possible forms of the leapfrog: the drift-kick-drift, also known as position-Verlet, and kick-drift-kick, known as velocity-Verlet,\footnote{\change{Our nomenclature, [DK]$^2$ for the drift-kick-drift leapfrog \change{describes} its composition as \change{symplectic Euler drift-kick (=DK) for $h/2$ followed by its adjoint for another $h/2$.} We use this scheme to name all maps in this study.}}
\begin{subequations} \label{eq:map:leapfrog}
\begin{eqnarray}
	\label{eq:map:[DK]^2}
	\change{\map{h}{H} \to\;} \psi_h^{\mathrm{\change{[DK]^2}}} &\equiv&
		\Exp{\frac{h}{2}\op{T}}\,\map{h}{V}\,\Exp{\frac{h}{2}\op{T}},
	\\ \label{eq:map:[KD]^2}
	\change{\map{h}{H} \to\;} \psi_h^{\mathrm{\change{[KD]^2}}} &\equiv&
		\Exp{\frac{h}{2}\op{V}}\,\map{h}{T}\,\Exp{\frac{h}{2}\op{V}}
\end{eqnarray}
\end{subequations}
with error Hamiltonians
\begin{subequations} \label{eq:leapfrog:H:err}
\begin{eqnarray}
	\label{eq:map:[DK]^2:H:err}
	\sub{H}{err}^{\mathrm{\change{[DK]^2}}} &=&
	- \frac{h^2}{24}\{\opn V,T\cls,T\}
	+ \frac{h^2}{12}\{\opn T,V\cls,V\} + \mathcal{O}(h^4),
	\\[1ex] \label{eq:map:[KD]^2:H:err}
	\sub{H}{err}^{\mathrm{\change{[KD]^2}}} &=&
	\phantom{-} \frac{h^2}{12}\{\opn V,T\cls,T\}
	- \frac{h^2}{24}\{\opn T,V\cls,V\} + \mathcal{O}(h^4).
\end{eqnarray}
\end{subequations}

\subsection{The error terms of the Leapfrog} \label{sec:leapfrog:err}
Let us take a closer look at the $\mathcal{O}(h^2)$ error terms of the leapfrog integrator in equations~(\ref{eq:leapfrog:H:err}). As the potential energy $V$ is the sum of the contributions $V_{\ij}$ from each pair-wise interaction, so is the error term $\{\opn V,T\cls,T\}$ the sum over the terms
\begin{equation}
	\label{eq:VijTT}
	\{\opn V_{\ij},T\cls,T\} =  \frac{Gm_im_j}{r_{\ij}^5}
		\left[v_{\ij}^2\,r_{\ij}^2-3(\vec{v}_{\ij}\cdot\vec{x}_{\ij})^2\right],
\end{equation}
where $\vec{v}_{\ij}\equiv\vec{v}_i-\vec{v}_{\!j}$ is the velocity difference, while
$v_{\ij}\equiv|\vec{v}_{\ij}|$ and $r_{\ij}\equiv|\vec{x}_{\ij}|$. The terms~(\ref{eq:VijTT}) become large only in a close encounter between particles $i$ and $j$. Assuming a parabolic encounter, we have $\frac{1}{2}v_{\ij}^2=G(m_i+m_j)/r_{\ij}$ such that $\{\opn V_{\ij},T\cls,T\}$ has magnitude $\sim G^2m_im_j(m_i+m_j)/r_{\ij}^4$.

The second contribution to the errors in equations~(\ref{eq:leapfrog:H:err}) are sums over terms of the form $\{\opn T,V_{\ij}\cls,V_{lk}\}$. These vanish if all four indices differ, and the only non-vanishing terms are of two types:
\begin{subequations} \label{eqs:TVV}
\begin{eqnarray}
	\label{eq:TVijVij}
	\{\opn T,V_{\ij}\cls,V_{\ij}\} &=&
		\frac{G^2m_im_{\!j}(m_i+m_{\!j})}{r_{\ij}^4}
	\qquad\text{and}
	\\[1ex]
	\label{eq:TVijVik}
	\{\opn T,V_{\ij}\cls,V_{ik}\} &=&
		\frac{G^2m_im_{\!j}m_k}{r_{\ij}^3\,r_{ik}^3}
	\vec{x}_{\ij}\cdot\vec{x}_{ik}
	\qquad\text{with}\quad j\neq k.
\end{eqnarray}
\end{subequations}
The first of these becomes large only in close encounters between particles $i$ and $j$, when it is of the same magnitude as the term $\{\opn V_{\ij},T\cls,T\}$ above. The second type of terms (equation~\ref{eq:TVijVik}) becomes large only in a close three-body encounter between particles $i$, $j$, and $k$ (close encounters of more than three particles only contribute to \change{yet} higher-order error terms, see Appendix~\ref{app:Err:4}). In order to distinguish these different terms, we define
\begin{subequations} \label{eqs:TVV:2,3}
\begin{eqnarray}
	\{\opn T,V\cls,V\}_2 &\equiv& \sum_{{i<j}}
		\{\opn T,V_{\ij}\cls,V_{\ij}\}
\end{eqnarray}
and
\begin{eqnarray}
	\{\opn T,V\cls,V\}_3 &\equiv&
	\change{2\!\sum_{i<j<k}
		\{\{T,V_{\!\ij}\},V_{\!ik}\} +
		\{\{T,V_{\!jk}\},V_{\!\!ji}\} +
		\{\{T,V_{\!ki}\},V_{\!k\!j}\},
	}
\end{eqnarray}
\end{subequations}
such that $\{\opn T,V\cls,V\}=\{\opn T,V\cls,V\}_2+\{\opn T,V\cls,V\}_3$
\change{(see equations~\ref{eq:app:TVV} and \ref{eq:app:TVV3} for computationally more useful alternative expressions)}.

\subsection{Higher order symplectic integrators} \label{sec:higher}
It is well known that in order to construct higher than second-order integrators by operator splitting, i.e.\ by alternating kicks and drifts with step sizes chosen such that the $\mathcal{O}(h^2)$ terms are eliminated from $\sub{H}{err}$, one must perform at least on backward kick and one backward drift \citep{Sheng1989,Suzuki1991}. Such methods have been proposed \citep[e.g.][]{Yoshida1990} but are rarely used in astrophysics, because backward steps are problematic \change{with frictional forces (such as tidal dissipation)}, but also because the coefficients of the errors terms \change{tend to be quite} large.

However, in order to obtain a fourth-order method not both of the error terms \change{in equations~(\ref{eq:leapfrog:H:err})} need to be eliminated: the \change{second} of these
\begin{equation} \label{eq:G}
	G \equiv \{\opn T,V\cls,V\} = \sum_k
	\frac{1}{m_k} \pdiff{V}{\vec{x}_k}\cdot\pdiff{V}{\vec{x}_k}
\end{equation}
depends only on the positions and can be integrated (see also Appendix~\ref{app:G})\change{. I}n other words\change{,} the map $\map{h}{G}$ is exactly soluble. This allows the construction of fourth-order symplectic integrators with only forward steps \citep{Suzuki1995,Chin1997}. The simplest such integrator is based on the relation 
\begin{equation} \label{eq:BCH:4th}
	\log\left(\Exp{\frac{1}{6}X}\,\Exp{\frac{1}{2}Y}\,\Exp{\frac{2}{3}X}\,
		 \Exp{\frac{1}{2}Y}\,\Exp{\frac{1}{6}X}\right)
	= X+Y +\tfrac{1}{72}[X,\change{[}X,Y\change{]}] + \dots,
\end{equation}
which implies that the map
\begin{equation} \label{eq:map:CC:2}
	\change{
		\psi_h^{\mathrm{[KDK]^2}} \equiv \;
	}
	\Exp{\frac{ h}{6}\op{V}}\,
	\Exp{\frac{ h}{2}\op{T}}\,
	\Exp{\frac{2h}{3}\op{V}}\,
	\Exp{\frac{ h}{2}\op{T}}\,
	\Exp{\frac{ h}{6}\op{V}}
\end{equation}
has error Hamiltonian
\begin{equation} \label{eq:CC:2:H:err}
	\sub{H}{err} = \frac{h^2}{72}\{\opn T,V\cls,V\}
		+ \mathcal{O}(h^4).
\end{equation}
Combining \change{(\ref{eq:map:CC:2})} with the map $\Exp{-\frac{h^3}{72}\op{G}}$ obtains the fourth-order forward integrator (dubbed `4A' by \citealt{Chin1997}, see also \citealt{ChinChen2005})
\begin{equation} \label{eq:CC:4}
	\map{h}{H} \to \psi_h^{\change{\mathrm{[KDK]^2_4}}} \equiv
	\Exp{\frac{ h}{6}\op{V}}\,
	\Exp{\frac{ h}{2}\op{T}}\,
	\Exp{\frac{2h}{3}(\op{V}-\frac{h^2}{48}\op{G})}\,
	\Exp{\frac{ h}{2}\op{T}}\,
	\Exp{\frac{ h}{6}\op{V}}.
\end{equation}
Here, the map $\Exp{-\frac{h^3}{72}\op{G}}$, corresponding to a force-gradient kick, is applied in the middle, such that the integrator remains self-adjoint, but for the order of the method this does not matter as long as it is applied at any time during the step.

More general symplectic maps can be constructed by alternating application of drifts, kicks, and force-gradient kicks. By carefully arranging the sub-steps of these component maps, the coefficients of the $\mathcal{O}(h^4)$ contributions to $\sub{H}{err}$ can be substantially reduced compared to the map~(\ref{eq:CC:4}) \citep*{OmelyanMryglodFolk2002, OmelyanMryglodFolk2003}. However, in order to obtain a sixth-order integrator, i.e.\ to have vanishing coefficients for all the $\mathcal{O}(h^4)$ contributions to $\sub{H}{err}$, backward steps are required, unless the error term $\{\change{\{\{\{}V,T\cls,T\cls,T\cls,V\}$ can be integrated \citep{Chin2005}, which is generally impossible.

\section{Symplectic maps using a Kepler solver} \label{sec:kepler}
Recently \cite*{GoncalvesFerrariEtAl2014} proposed to replace the pair-wise kick map $\Exp{h\op{V}_{\ij}}$ for each particle pair with \change{a backwards drift followed by their mutual binary orbit}, hereafter a \emph{binary kick}:
\begin{equation} \label{eq:binary:kick}
	\Exp{h\op{V}_{\ij}} \to
	\Exp{h\op{H}_{\ij}}\,\Exp{-h(\op{T}_i+\op{T}_{\!j})}
\end{equation}
\citep{HernandezBertschinger2015} with the binary Hamiltonian
\begin{equation} \label{eq:H:binary}
	H_{\ij} \equiv T_i + T_{\!\!j} + V_{\!\ij}
	 = \frac{\vec{p}_{\!i}^2}{2m_{\!i}} + \frac{\vec{p}_{\!j}^2}{2m_{\!j}} -
		\frac{Gm_im_{\!j}}{r_{\ij}}.
\end{equation}
Since the forward drift of the centre of mass due to $\Exp{h\op{H}_{\ij}}$ cancels with its backward drift due to $\Exp{-h(\op{T}_{i}+\op{T}_{\!j})}$, the map~(\ref{eq:binary:kick}) can also be implemented via the equivalent form
(\citeauthor{GoncalvesFerrariEtAl2014})
\begin{equation} \label{eq:binary:kick:alt}
	\Exp{h\op{K}_{\ij}}\, \Exp{-h\{.,\frac{1}{2}\mu_{\ij}^{}\vec{v}_{\ij}^2\}},
\end{equation}
where $\mu_{\ij}\equiv m_im_j/(m_i+m_j)$ is the reduced mass, and
\begin{equation}
	K_{\ij} \equiv
	\mu_{\ij} \left[\frac{v_{\ij}^2}{2} - \frac{G(m_i+m_j)}{r_{\ij}}
	\right]
\end{equation}
the Kepler Hamiltonian of the particle pair. \change{We} found no detectable difference between maps \eqref{eq:binary:kick:alt} and \eqref{eq:binary:kick} in terms of computational efficiency or finite precision errors.

\subsection{The method of Hernandez \& Bertschinger revisited} \label{sec:HB15}
The map defined in equation~(\ref{eq:binary:kick}) or~(\ref{eq:binary:kick:alt}) is not self-adjoint, such that substituting it for every pair-wise kick\change{, i.e.\ replacing} 
\begin{equation} \label{eq:map:W}
	\change{
		\map{h}{V} \to\;\;
	}
	\psi_{h}^{W}
	\equiv \prod_{\change{\text{$(i,j)$ in some order}}}
		\Exp{h\op{H}_{\ij}}\,\Exp{-h(\op{T}_i+\op{T}_{\!j})},
\end{equation}
in the standard $N$-body integrators~(\ref{eq:map:leapfrog}) obtains a method that is not self-adjoint either and hence also not \change{reversible} (\citeauthor{HernandezBertschinger2015}). A self-adjoint integrator can be composed \change{from any map $\altmap{h}{}$ as $\psi_h=\altmap{h/2}{\dag}\altmap{h/2}{}$. \citeauthor{HernandezBertschinger2015} applied this recipe to the irreversible map \begin{equation} \label{eq:map:phi}
	\phi_h \equiv \psi_{h}^{W} \map{h}{T},
\end{equation}
which is similar to the symplectic Euler but second-order accurate (\citeauthor{HernandezBertschinger2015}). This yields the integrator}
\begin{subequations} \label{eq:map:HB15}
\begin{equation}\label{eq:map:[DB]^2}
	\map{h}{H} \to \psi_h^{\mathrm{\change{[DB]^2}}} \equiv
	\change{
		\phi_{h/2}^{\dag}\phi_{h/2}^{\phantom{\dag}} = \;\;
	}
	\Exp{\frac{h}{2}\op{T}}\,\psi_{h/2}^{\dag W}\,\psi_{h/2}^W\,\Exp{\frac{h}{2}\op{T}},
\end{equation}
hereafter `HB15' \change{or [DB]$^2$} \change{with `B' for binary kick.}
Alternatively, \change{the reversed recipe puts} the drift operation in the middle:
\begin{equation}\label{eq:map:[BD]^2}
	\map{h}{H} \to
	\psi_h^{\mathrm{\change{[BD]^2}}} \equiv
	\change{
		\phi_{h/2}^{\phantom{\dag}}\phi_{h/2}^{\dag} = \;\;
	}
	\psi_{h/2}^{W}\,\map{h}{T}\,\psi_{h/2}^{\dag W}.
\end{equation}
\end{subequations}
These integrators look different and would definitely be implemented differently\change{. However, as shown in equation~\eqref{eq:map:[DB]:recursive}, the maps $\phi_h^{}$ and $\phi_h^{\dag}$ are identical save for a reversal of the order of binary kicks. Hence, the maps~\eqref{eq:map:[DB]^2} and \eqref{eq:map:[BD]^2} differ only by a reversal of the binary-kick order in each half.} \citeauthor{HernandezBertschinger2015} reported that this order has no significant effect on the accuracy of the method, and indeed it does not affect the error Hamiltonian at order $\mathcal{O}(h^2)$. In Appendix~\ref{app:H:err:HB15}, we derive the error Hamiltonian for the maps~(\ref{eq:map:[DB]^2}) and~(\ref{eq:map:[BD]^2}) to be
\begin{equation} \label{eq:BDB:H:err}
	\sub{H}{err}^{\mathrm{\change{[DB]^2}}}
	= \frac{h^2}{48}\{\opn T,V\cls,V\}_3
	+ \mathcal{O}(h^4)
\end{equation}
and $\sub{H}{err}^{\mathrm{\change{[BD]^2}}}=\sub{H}{err}^{\mathrm{\change{[DB]^2}}}+\mathcal{O}(h^4)$. Comparing this to the errors of the leapfrog, as discussed in section~\ref{sec:leapfrog:err}, we see that all $\mathcal{O}(h^2)$ error terms arising from close two-body encounters have been removed. The only remaining terms are of the form~(\ref{eq:TVijVik}), which account for close three-body encounters. In fact, this is true at all orders: error terms which are (nested) Poisson brackets containing only $T$ and $V_{\ij}$ are eliminated at all orders. This can be seen by letting $m_k\to0$ for all but one pair of particles, when the method becomes exact for this pair while the Hamiltonian collapses to~(\ref{eq:H:binary}). Thus, the only remaining error terms \change{involve two or more particle pairs}.

A more general self-adjoint arrangement of the maps~(\ref{eq:map:W}) and $\map{h}{T}$ is
\begin{equation} \label{eq:map:psi:alpha}
	\map{h}{H} \to \psi^\alpha_h \equiv
	\map{\alpha h}{T}\,\psi_{h/2}^{W}\,\map{(1-2\alpha)h}{T}\,
		\psi_{h/2}^{\dag W}\,\map{\alpha h}{T}
\end{equation}
with parameter $\alpha$. For $\alpha=0$, we obtain the map \change{[BD]$^2$}, while $\alpha=\tfrac{1}{4}$ corresponds to the integrator $\zeta_h^2$ of \citeauthor{HernandezBertschinger2015}. For $N=2$, \change{map~\eqref{eq:map:psi:alpha}} reduces to a simple Kepler solver only \change{for} $\alpha=0$ \change{and,} consequently, only this choice eliminates both the error terms $\{\opn V,T\cls,T\}$ and $\{\opn T,V\cls,V\}_2$ (we give $\sub{H}{err}$ up to order $h^2$ in equation~\ref{eq:err:map:psi:alpha}). \change{This} explains why \citeauthor{HernandezBertschinger2015} found their $\zeta_h^2$ integrator to be inferior to \change{[DB]$^2$}.

\subsection{Extending the method of Hernandez \& Bertschinger to fourth order} \label{sec:HB15:4}
The error Hamiltonian~(\ref{eq:BDB:H:err}) to second order is in fact integrable and fourth-order schemes can be constructed by simply integrating it, analogously to the forward fourth-order integrators discussed in section~\ref{sec:higher}. The costs for doing so are small compared to those for the solutions to the $N(N-1)/2$ Kepler problems, though both scale as $\mathcal{O}(N^2)$. The fourth-order correction can be applied either in the middle, beginning and end, or both. Describing this freedom with parameter $\alpha$, gives the fourth-order map
\begin{equation}
	\label{eq:map:[DB]^2_4}
	\psi_h^{\mathrm{\change{[DB]^2_4}}} \equiv
	\Exp{-\alpha\frac{h^3}{96}\op{G}_s}\;
	\change{
		\phi_{h/2}^{\dag}\;
	}
	\Exp{(\alpha-1)\frac{h^3}{48}\op{G}_s}\;
	\change{
		\phi_{h/2}^{\phantom{\dag}}\;
	}
	\Exp{-\alpha\frac{h^3}{96}\op{G}_s},
\end{equation}
where $G_s \equiv \{\opn T,V\cls,V\}_3$.
\change{Again, swapping the sub-steps $\phi_{h/2}^{\dag}$ and $\phi_{h/2}^{}$ obtains a map, [BD]$^2_4$, which is identical except for a reversal of the binary-kick order in each half.}
\change{The map~\eqref{eq:map:[DB]^2_4} and its generalisation in equation~\eqref{eq:map:[KDBK]^2_4} below are the main results of this study and we will also call them `DH16'.}

\begin{figure}
	\begin{center}
	\includegraphics[width=0.85\columnwidth]{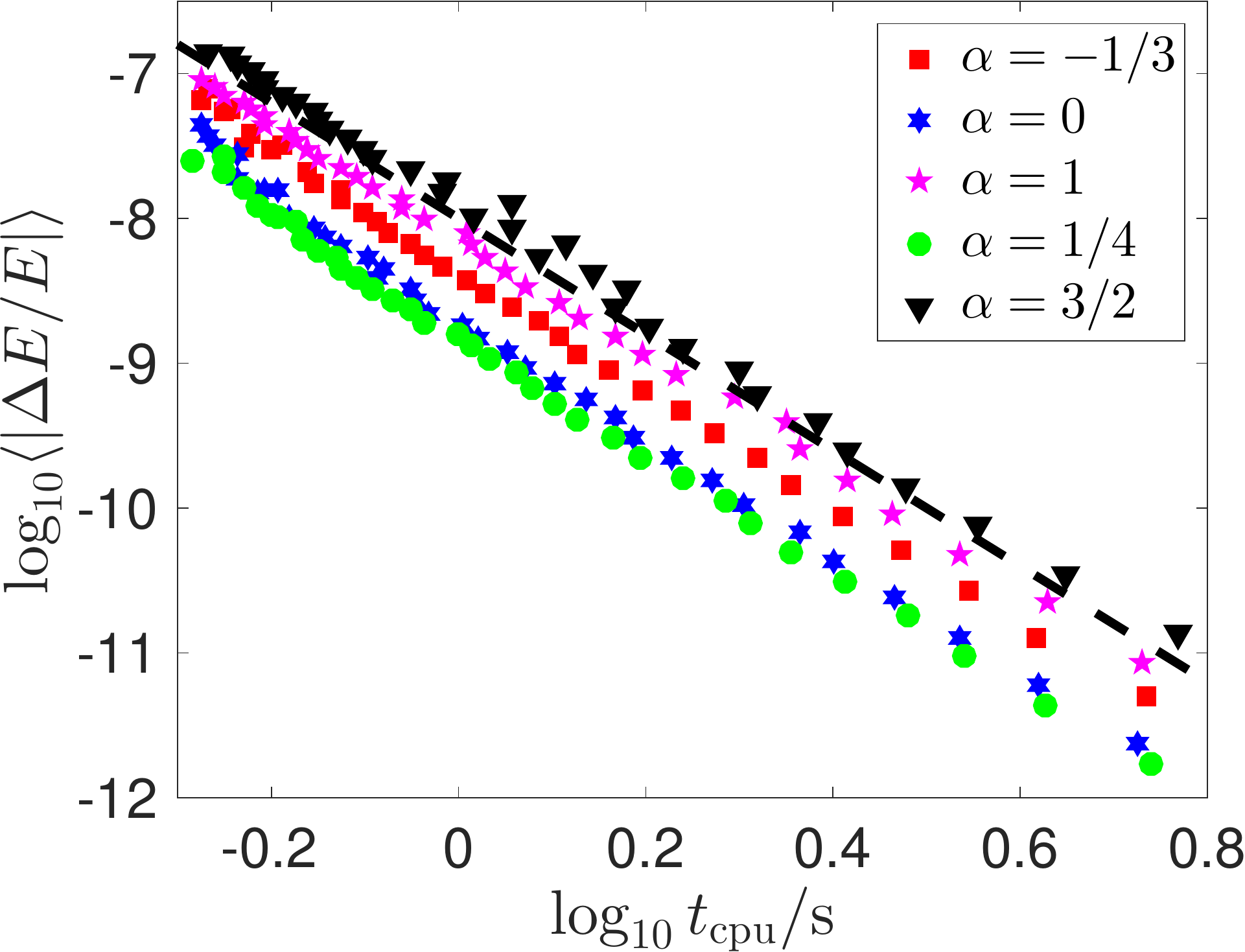}
	\end{center}
	\caption{Mean absolute energy error plotted vs.\ computational costs for an integration of the outer Solar system (Sun and gas giants) over 1000 years using the map~\eqref{eq:map:[DB]^2_4} for various choices of time step $h$ and parameter $\alpha$. The line indicates $\sub{t}{cpu}^{-4}$, the scaling expected for a fourth-order method.
	\label{fig:alpha}
	}
\end{figure}

To confirm experimentally that the map~\eqref{eq:map:[DB]^2_4} is indeed fourth-order accurate, we integrate the Sun and outer gas giant planets \cite[with initial conditions taken from][]{HairerLubichWanner2006} for $t=1000$ years. For given step size $h$, we iterate the map~\eqref{eq:map:[DB]^2_4} $n=\lfloor t/h\rfloor$ times, and calculate the arithmetic mean $\langle|\Delta E/E|\rangle$ of the energy error magnitudes after each step. We also measure the total computational time $\sub{t}{cpu}$, which we expect to scale as $\sub{t}{cpu}\sim 1/h$ \citep[but see ][for \change{a} discussion on variations in the cost of the Kepler solver depending on $h$]{WisdomHernandez2015} and repeat these calculations for various values of $h$ and $\alpha$. For any given $\alpha$, we find in Fig.~\ref{fig:alpha} that $\langle|\Delta E/E|\rangle\propto\sub{t}{cpu}^{-4}$ as expected for a fourth-order method.

Fig.~\ref{fig:alpha} also shows that $\alpha=1/4$ performs better than $\alpha=0$, despite the extra force-gradient map (the Kepler solver dominates the computational costs). In Appendix \ref{sec:fourtherror}, we derive $\sub{H}{err}$ for the integrator~\eqref{eq:map:[DB]^2_4} to fourth order in $h$. While at this order $\sub{H}{err}$ is a linear function of $\alpha$, the various error terms depend non-trivially on the state of the $N$-body system as well as the ordering of particle pairs within the binary-kick operator $\psi_h^W$. This makes it very difficult, if not impossible, to deduce the optimal $\alpha$ and ordering of pairs by analysis alone. Instead, we explore \change{the} optimal $\alpha$ numerically.

\begin{figure}
	\begin{center}
	\includegraphics[width=0.85\columnwidth]{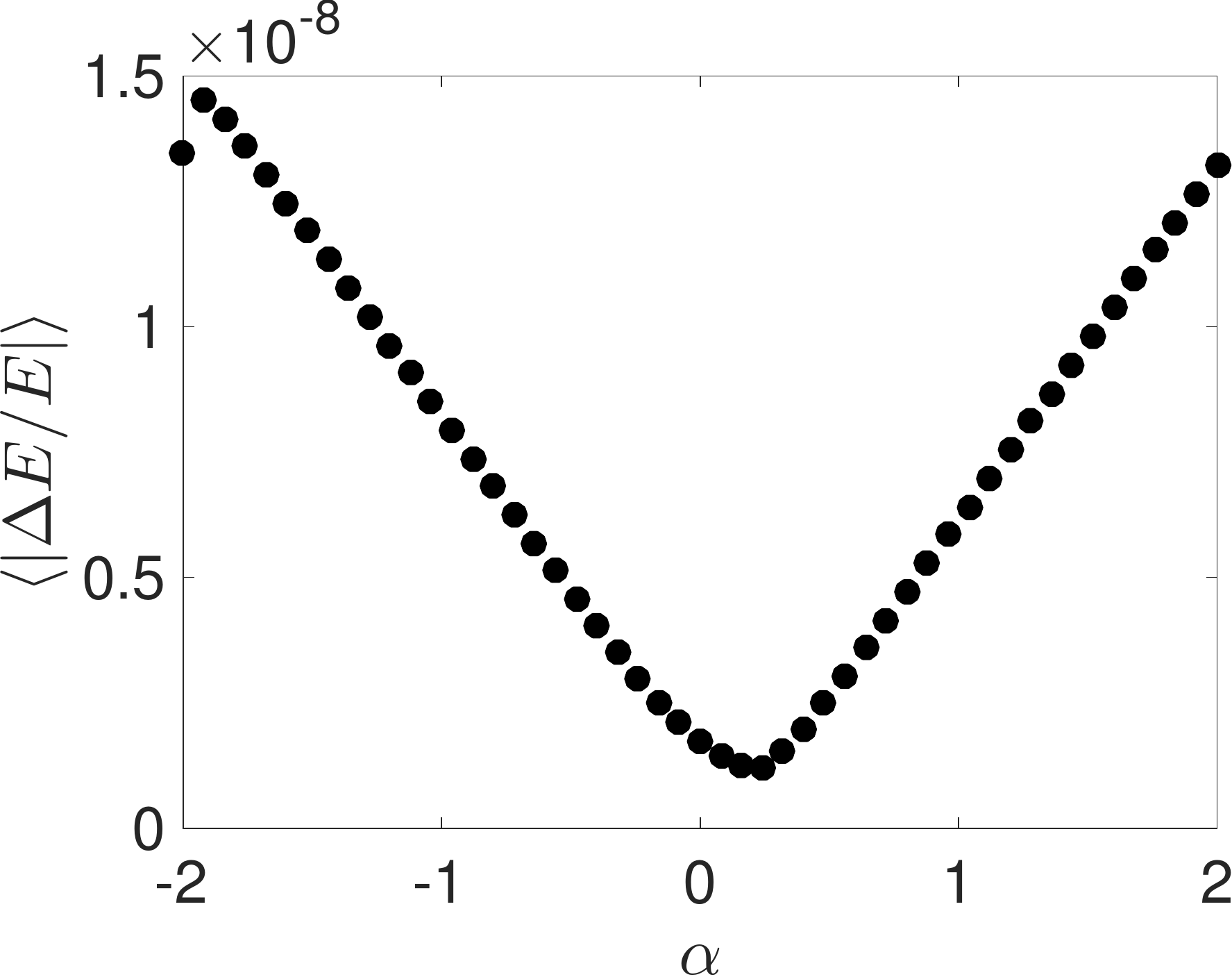}
	\end{center}
	\caption{Mean absolute energy error plotted vs.\ parameter $\alpha$ for integrations as shown in Fig.~\ref{fig:alpha} but at fixed time step $h=0.5\,$yr. For this experiment, the particle pair ordering in the binary-kick operator $\psi_h^W$ was fixed at the Sun-planet pairs followed by the planet-planet pairs The optimal $\alpha$ depends slightly on the pair ordering and the $N$-body problem itself.
	\label{fig:alpha30}
	}
\end{figure}

To this end, we repeat the previous calculations at fixed $h=0.5\,$yr, when we find
$\sub{t}{cpu}$ to hardly depend on $\alpha$. Fig.~\ref{fig:alpha30} plots the mean energy error magnitude versus $\alpha$. For our test problem of the outer Solar system, we find that $\sub{\alpha}{opt}=1/4$ obtains the smallest error. For the $\alpha$ ranges in this plot, $\langle|\Delta E/E|\rangle$ varies by about a factor 10. However, this does depend slightly on the ordering of the particle pairs within $\psi_h^W$. Varying the order of the Sun-planet pairs or that of the planet-planet pairs leaves $\sub{\alpha}{opt}$ approximately unchanged. However, when we reversed fully the pair ordering, $\sub{\alpha}{opt}\approx0.16$.

We also measured $\sub{\alpha}{opt}$ for different $N$-body problems. For the hierarchical triple problem considered by \cite{DuncanLevisonLee1998} and \cite{Hernandez2016} we find still $\sub{\alpha}{opt} \approx 1/4$ for $h=0.001\,$yr and $t$=1\,yr. For the figure-of-eight three-body solution discussed in \cite{ChencinerMontgomery2000}, we find $\sub{\alpha}{opt}\approx 0.15$ for $t$ equal to the period P and $h=P/50$. In summary, the optimal value of $\alpha$ for the integrator \eqref{eq:map:[DB]^2_4} appears to vary depending on the $N$-body problem and the solution strategy but we always find it to be constrained between 0.1 and 0.3.
  
\subsection{A test of a chaotic exchange orbit}
\label{sec:chaotex}
We now test our fourth-order map~\eqref{eq:map:[DB]^2_4} on a challenging problem: a chaotic exchange orbit of the planar restricted circular three-body problem. If we denote the coordinates and velocities in the co-rotating frame with a prime, then the Jacobi integral (the only isolating integral for this problem) is
\begin{equation}
	C_J
	= \tfrac{1}{2}\vec{v}^2 + \Phi(\vec{x}') - \vec{\omega}\cdot(\vec{x}\times\vec{v})
	= \tfrac{1}{2}\vec{v}^{\prime2} + U(\vec{x}')
\end{equation}
with binary angular velocity $\vec{\omega}=\sqrt{G(m_1+m_2)/a^3}\hat{\vec{z}}$,
\begin{equation}
	\Phi(\vec{x}') = -\frac{Gm_1}{|\vec{x}'-\vec{r}_1|}
					-\frac{Gm_2}{|\vec{x}'-\vec{r}_2|},
\end{equation}
and $U(\vec{x}')\equiv\Phi(\vec{x}')-\tfrac{1}{2}(\vec{\omega}\times\vec{x}')^2$. The conventional Jacobi constant definition is $C = - 2 C_J$, but $C_J$ is equal in value to the Hamiltonian in the co-rotating frame. Here, $a$ is the binary semi-major axis, $m_{1,2}$ the masses of its components, and $\vec{r}_{1,2}$ their co-rotating positions. If $\vec{L}_{1,2}$ \change{are} the co-rotating positions of the first and second Lagrange points, then orbits satisfying $U(\vec{L}_1)<C_J <U(\vec{L}_2)$ can visit both masses but cannot escape to infinity.

We use units of au, days, and Solar mass, when we set $m_1=1$, $\mu = m_2/(m_1+m_2)=0.01$, and $a=5.2$. We integrate the orbits of all three particles in the barycentric inertial frame, starting the binary components on the $x$-axis and the test particle at (4.42,0,0) and with velocity (0,0.0072,0), both w.r.t.\ the Solar mass object.  With respect to the center of mass the position coordinates are $\approx (4.369,0,0)$ and the velocity coordinates are $\approx (0,0.0071,0)$.  For these settings $C_J = -9.0770\times 10^{-5}$, $U(\vec{L}_1) = -9.1038 \times 10^{-5}$, and $U(\vec{L}_2) = -9.0654 \times 10^{-5}$, satisfying the conditions for a chaotic exchange orbit. 

\begin{figure}
	\begin{center}
	\includegraphics[width=0.85\columnwidth]{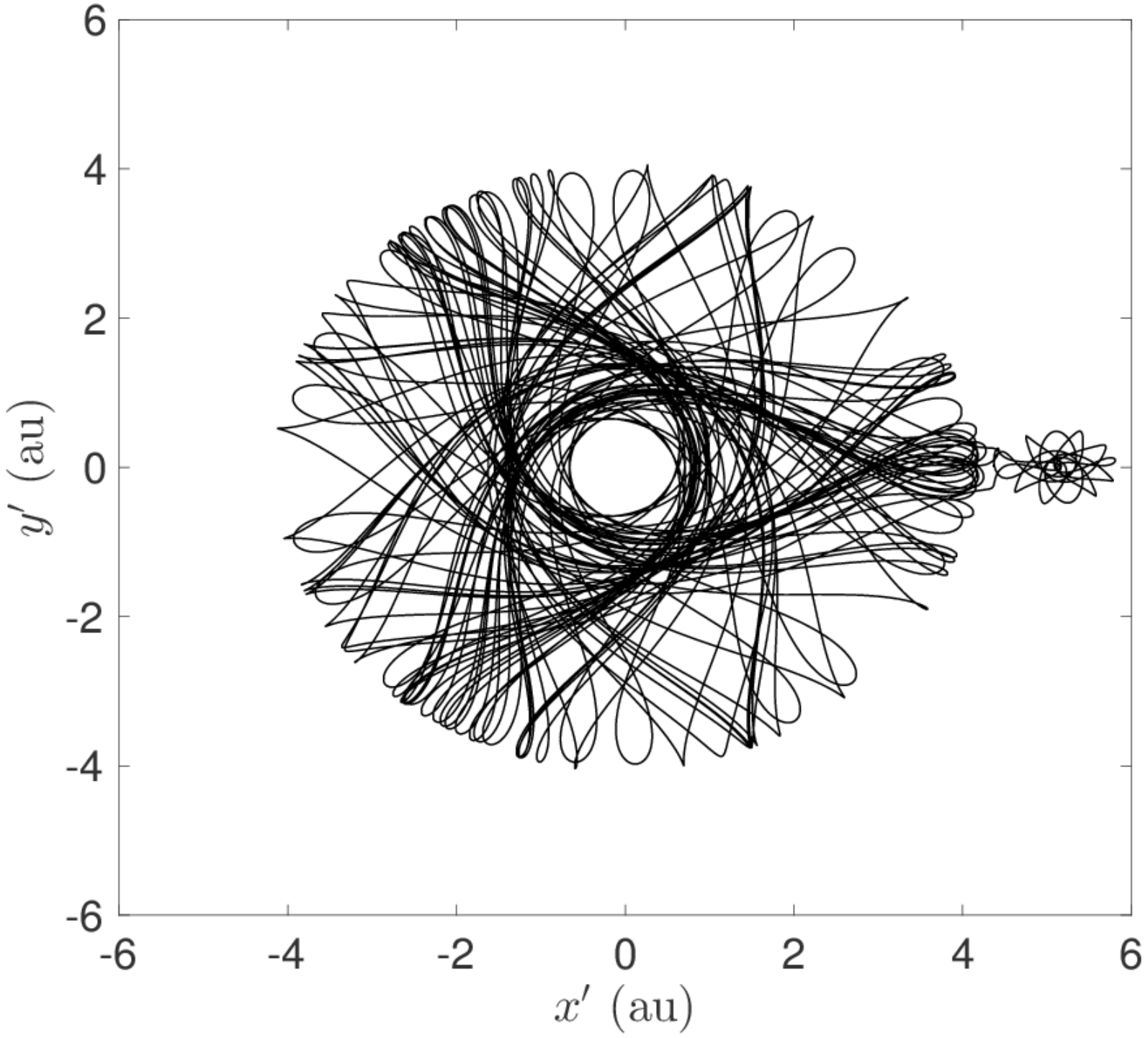}
	\end{center}
	\caption{Trajectory over 500 years of the test particle in the circular restricted three-body problem considered in the text. Most of the time the test particle orbits the primary, but occasionally switches to the secondary.
	\label{fig:exchangeorbit}
  	}
\end{figure}
\begin{figure}
	\begin{center}
	\includegraphics[width=\columnwidth]{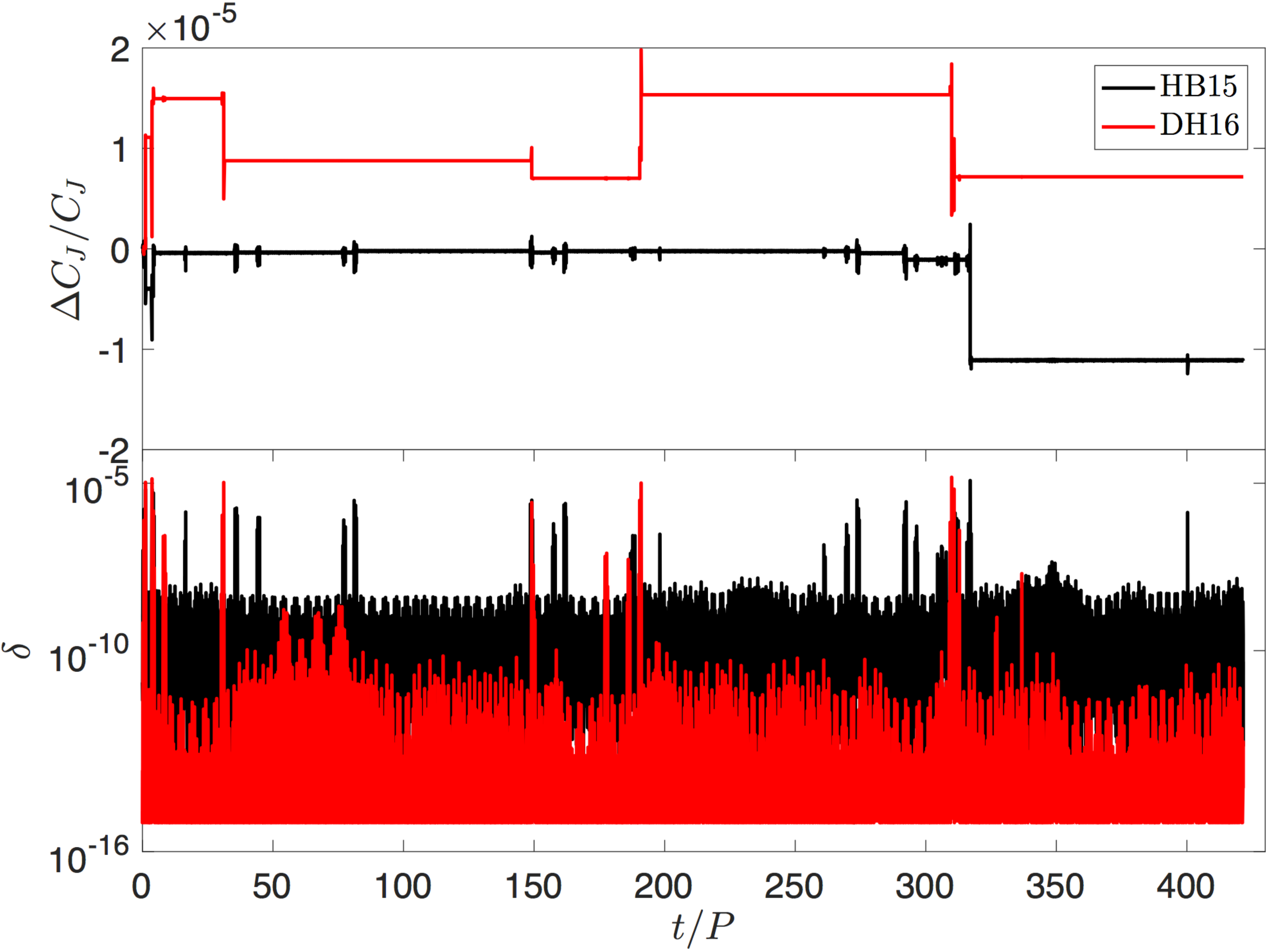}
	\end{center}
	\caption{\change{Relative e}rror in the Jacobi constant effected by \change{[DB]$^2$ (}HB15\change{)} and \change{[DB]$^2_4$ with $\alpha=1$ (}DH16\change{)}. \change{The top panel gives the accumulated error, while in the bottom panel the error $\delta=|C_J(t)-C_J(t-h)|/|C_J(t-h)|$ over one time step is plotted.} The initial conditions are the chaotic exchange orbit described in the text. \change{While the local error for the fourth-order method is considerably better, the accumulated} error magnitudes of the two \change{maps} are generally similar. Due to the chaotic nature of the orbit, the actual trajectories \change{of the two integrations} differ after $t\approx5\,P$.
	\label{fig:restricted}
	}
\end{figure} 
Fig.~\ref{fig:exchangeorbit} plots the trajectory of the test particle over 500 years in the co-rotating frame. The test particle has multiple close encounters with $m_2$ within its Hill radius.

The period of the massive bodies is $P\approx11.9$ years.  We compute the Lyapunov time for this problem using map \eqref{eq:map:[DB]^2_4} since it is the map we are interested in studying.  But note that the Lyapunov time can be a function of the map and $h$.  For $h = 0.1$ y\change{ea}r\change{s} and nearby initial conditions we calculate a Lyapunov time $t_{\mathrm{L}}\approx0.3\,P$. Fig.~\ref{fig:restricted} plots the error of the Jacobi integral as a function of time for an integration over $5000\text{ years}=421\,P=1405\,t_{\mathrm{L}}$ using $h=4$ days.

We see that the \change{accumulated} errors of \change{[DB]$^2$} and \change{[DB]$^2_4$} are similar, \change{though the local error of the fourth-order method is substantially smaller, often reaching the round-off limit of $\sim10^{-15}$}. We have found in other experiments that \change{in the presence of close encounters the performance of symplectic integrators can deteriorate}, and it is not surprising \change{[DB]$^2_4$} and \change{[DB]$^2$} behave similarly.

If we let $h=8\,$days, a specialised integrator for Solar system problems, \texttt{MERCURY} \citep{Chambers1999}, yields an error of the order $10^{-5}$. \texttt{MERCURY} has been found to not always be symplectic and tends to yield wrong behaviour for three-body problems \citep{Hernandez2016}. For a contrasting example, consider the forward stepping fourth-order map~(\ref{eq:CC:4}), corresponding to map~\eqref{eq:map:[KDBK]^2_4} with $\mathcal{S}$ empty, i.e.\ without employing a Kepler solver. If we set $h=4$ and a short $t=20$ years, this map yields a large error $|\Delta C_J/C_J| = 0.047$.

If we let \change{[DB]$^2_4$} run longer than 5000 years (but still at $h=4\,$days), the Jacobi energy error may jump by orders magnitude, whereas \texttt{MERCURY} does not yield such jumps. The jumps are associated with close encounters to $m_2$. However, while some close approaches caused jumps, other closer approaches did not. This indicates that a constant time step of $h=4\,$days is inappropriate for this problem in the long term.

\change{
\subsection{An \boldmath $N$-body test}\label{sec:nbody}
We also test the maps~\eqref{eq:map:[DB]^2} and \eqref{eq:map:[DB]^2_4} with $\alpha=0$ for an $N$-body system. To this end, we use an implementation dubbed \textsc{triton} which employs computational parallelism (see appendix~\ref{app:eff:Kepler} for details). We simulate a cluster of $N=1024$ equal-mass particles, initially following a \cite{Plummer1911} model with ergodic distribution function, equivalent to the simulations reported in Fig.~5 of \citeauthor{GoncalvesFerrariEtAl2014}. Like those authors, we use $N$-body units ($G=1$, $M=1$, and $E=-1/4$, which imply a virial radius of 1 and a crossing and relaxation time at half-mass of $\sim2.4$ and $\sim45$, respectively) and integrate the system from $t=0$ to $t=400$ with steps of $h=10^{-4}$. Our initial conditions are different from those used by \citeauthor{GoncalvesFerrariEtAl2014}, but equivalent in the sense that we use the same model to draw them from (we set the centre of mass and total momentum to zero).

\begin{figure}
	\begin{center}
	\includegraphics[width=\columnwidth]{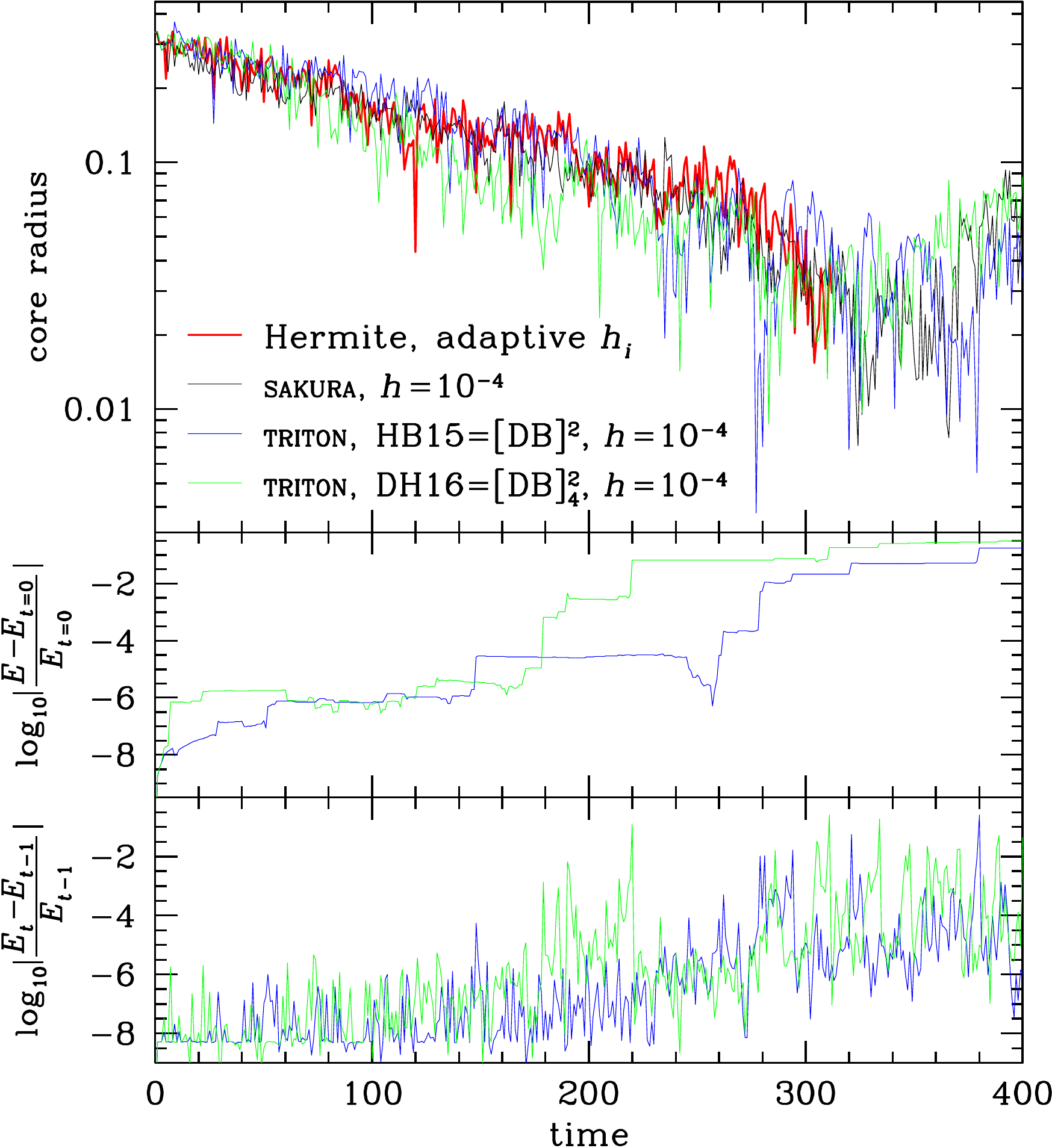}
	\end{center}
	\caption{\change{Core radius \citep*[calculated as proposed by][]{McMillanHutMakino1990} and energy errors for simulations of a 1024-body Plummer sphere. The core-radii data for the time-adaptive Hermite integrator and \textsc{sakura} \citep{GoncalvesFerrariEtAl2014} are taken from Fig.~5 of \citeauthor{GoncalvesFerrariEtAl2014} (the Hermite integrator ground to a hold at core collapse because of close encounters). The middle and bottom panels plot, respectively, the relative accumulated and short-term energy error for [DB]$^2$ (HB15) and [DB]$^2_4$ with $\alpha=0$ (DH16) only, which we implemented in parallel in a code \textsc{triton} as described in appendix~\ref{app:eff:Kepler}.}
	\label{fig:nbody}
	}
\end{figure} 
The top panel of Fig.~\ref{fig:nbody} plots the core radius as function of time for our runs as well as two simulations reported by \citeauthor{GoncalvesFerrariEtAl2014}: one with their code \textsc{sakura} also using $h=10^{-4}$ and another with a Hermite integrator using adaptive time stepping. \textsc{sakura}, which also uses a binary kicks, violates both symplecticity and time reversibility \citep{HernandezBertschinger2015}, but becomes exact in the two-particle limit, like the maps [DB]$^2$ and [DB]$^2_4$. Because of this, its truncation errors (which cannot be represented by an error Hamiltonian) are unlikely to contain contributions arising from two-body encounters. This property (which our maps share) enables a reasonably accurate integration through core collapse. There is no appreciable difference between the core-radius evolution of \textsc{sakura} and our maps.

The accumulated energy errors for our fixed-time step integration are considerable, reaching $>10\%$ at the final time, though staying at the same level of $|\delta E/E|\lesssim0.01$ until core collapse ($t\sim300$) as reported for \textsc{sakura}. There is no advantage of the fourth-order method. This is because the constant time step is simply too long to resolve close three-body encounters, which destroy any advantage of the fourth-order method and cause sudden increases of the accumulated error. The energy over a period of one time unit (bottom panel of Fig.~\ref{fig:nbody}) is much better behaved, though not surprisingly has increased by $\sim10^{2-3}$ by $t\sim300$, the time of core collapse.

We thus conclude from this test that \textsc{triton} is as good as \textsc{sakura} in its ability to integrate collision-dominated $N$-body dynamics, but unlike \textsc{sakura} is symplectic and reversible. Despite the necessity to make 
two calls to the kepler solver per particle pair and time step as opposed to 
\textsc{sakura}'s one, \textsc{triton} is about twice as fast (see appendix~\ref{app:eff:Kepler} for the likely reason).

}
\section{Using the Kepler solver for selected interactions only} \label{sec:mix}
The main problem with the method of the previous section is the computational expense of the Kepler solver needed in the binary kicks (but see Appendix~\ref{app:eff:Kepler}). The benefit from using such an approach is really only justified in close encounters. As already discussed by \citeauthor{HernandezBertschinger2015}, a faster method can be constructed by restricting the Kepler solver to selected pair-wise interactions. Let $\mathcal{S}$ be a set of $K\le N(N-1)/2$ particle pairs for which binary kicks shall be applied. Then we can split the potential energy into contributions integrated with binary kick and without:
\begin{equation}
	V = V_{\!s} + V_{\!c}
	\quad\text{with}\quad
	V_{\!s} = \sum_{(i,j)\in\mathcal{S}}V_{\!\ij}
	\quad\text{and}\quad
	V_{\!c} = \sum_{(i,j)\not{\in}\mathcal{S}}V_{\!\ij}.
\end{equation}
In analogy to equation\change{s}~\eqref{eq:map:W} \change{and \eqref{eq:map:phi}} we define the map\change{s}
\begin{equation} \label{eq:map:W:c}
	\psi_{h}^{W_s} \equiv \prod_{\text{$(i,j)\in\mathcal{S}$ \change{in some order}}}
	\Exp{h\op{H}_{\ij}}\,\Exp{-h(\op{T}_i+\op{T}_{\!j})}
\change{
	\qquad\text{and}\qquad
	\phi_h^s\equiv\psi_h^{W_s}\map{h}{T}
}
	.	
\end{equation}

\subsection{Extending the leapfrog}\label{sec:mix:2}
There are four di\change{stinc}t self-adjoint ways in which one can combine the maps $\Exp{h\op{V}_{\!c}}$, $\Exp{h\op{T}}$, and $\psi^{W_s}$ into a second-order integrator
\footnote{\change{Again,} alternatives \change{obtained by swapping $\phi_{h/2}^{\dag s}$ and $\phi_{h/2}^{s\phantom{\dag}}$ in equations~\eqref{eq:map:[DBK]^2} and \eqref{eq:map:[KDB]^2}} are identical \change{except for a reversal of the order of binary kicks}.}
\begin{subequations} \label{eq:maps:ex:leapfrog}
\begin{eqnarray}
	\label{eq:map:[DBK]^2}
	\psi_h^{\mathrm{\change{[\change{DB}K]^2}}} &=&
\change{
	\phi_{h/2}^{\dag s}\;
	\Exp{h\op{V}_{\!c}}\;
	\phi_{h/2}^{s\phantom{\dag}}
	=\;\;
}
	\psi_{h/2}^{W_s}\;\Exp{\frac{h}{2}\op{T}}\;
	\Exp{h\op{V}_{\!c}}\;
	\Exp{\frac{h}{2}\op{T}}\;\psi_{h/2}^{\dag W_s},
	\\
	\label{eq:map:[BKD]^2}
	\psi_h^{\mathrm{\change{[BKD]^2}}} &=&
	\psi_{h/2}^{W_s}\;\Exp{\frac{h}{2}\op{V}_{\!c}}\;\Exp{h\op{T}}\;
	\Exp{\frac{h}{2}\op{V}_{\!c}}\;\psi_{h/2}^{\dag W_s},
    \\
	\label{eq:map:[DKB]^2}
	\psi_h^{\mathrm{\change{[DKB]^2}}} &=&
	\Exp{\frac{h}{2}\op{T}}\;\Exp{\frac{h}{2}\op{V}_{\!c}}\;\psi_{h/2}^{\dag W_s}\;
	\psi_{h/2}^{W_s}\;\Exp{\frac{h}{2}\op{V}_{\!c}}\;\Exp{\frac{h}{2}\op{T}},
	\\
	\label{eq:map:[KDB]^2}
	\psi_h^{\mathrm{\change{[KDB]^2}}} &=&
\change{
	\Exp{\frac{h}{2}\op{V}_{\!c}}\;
	\phi_{h/2}^{\dag s}\;
	\phi_{h/2}^{s\phantom{\dag}}\;
	\Exp{\frac{h}{2}\op{V}_{\!c}}
	=\;\;
}
	\Exp{\frac{h}{2}\op{V}_{\!c}}\;\Exp{\frac{h}{2}\op{T}}\;\psi_{h/2}^{\dag W_s}\;
	\psi_{h/2}^{W_s}\;\Exp{\frac{h}{2}\op{T}}\;\Exp{\frac{h}{2}\op{V}_{\!c}},
\end{eqnarray}
\end{subequations}
with error Hamiltonians (derived in appendix~\ref{app:H:err:mix:2})
\begin{subequations} \label{eq:H:err:LF:ext}
\begin{eqnarray}
	\label{eq:H:err:[DBK]^2}
	\sub{H}{err}^{\mathrm{\change{[\change{DB}K]^2}}} &=&
		- \frac{h^2}{24}\{\opn V_{\!c},T\cls,T\}
		+ \frac{h^2}{12}\{\opn T,V_{\!c}\cls,V_{\!c}\}
		+ \frac{h^2}{24}\{\opn T,V_{\!s}\cls,V_{\!c}\}
		\nonumber \\ &&
		+ \frac{h^2}{48}\{\opn T,V_{\!s}\cls,V_{\!s}\}_3
		+ \mathcal{O}(h^4),
	\\[0.5ex]
	\label{eq:H:err:[BKD]^2}
	\sub{H}{err}^{\mathrm{\change{[BKD]^2}}} &=& \phantom{-}
		  \frac{h^2}{12}\{\opn V_{\!c},T\cls,T\}
		- \frac{h^2}{24}\{\opn T,V_{\!c}\cls,V_{\!c}\}
		+ \frac{h^2}{24}\{\opn T,V_{\!s}\cls,V_{\!c}\}
		\nonumber \\ &&
		+ \frac{h^2}{48}\{\opn T,V_{\!s}\cls,V_{\!s}\}_3
		+ \mathcal{O}(h^4),
	\\[0.5ex]
	\label{eq:H:err:[DKB]^2}
	\sub{H}{err}^{\mathrm{\change{[DKB]^2}}} &=&
		- \frac{h^2}{24}\{\opn V_{\!c},T\cls,T\}
		+ \frac{h^2}{12}\{\opn T,V_{\!c}\cls,V_{\!c}\}
		+ \frac{h^2}{6} \{\opn T,V_{\!s}\cls,V_{\!c}\}
		\nonumber \\ &&
		+ \frac{h^2}{48}\{\opn T,V_{\!s}\cls,V_{\!s}\}_3
		+ \mathcal{O}(h^4),
	\\[0.5ex]
	\label{eq:H:err:[KDB]^2}
	\sub{H}{err}^{\mathrm{\change{[KDB]^2}}} &=& \phantom{-}
		  \frac{h^2}{12}\{\opn V_{\!c},T\cls,T\}
		- \frac{h^2}{24}\{\opn T,V_{\!c}\cls,V_{\!c}\}
		- \frac{h^2}{12}\{\opn T,V_{\!s}\cls,V_{\!c}\}
		\nonumber \\ &&
		+ \frac{h^2}{48}\{\opn T,V_{\!s}\cls,V_{\!s}\}_3
		+ \mathcal{O}(h^4).
\end{eqnarray}
\end{subequations}
These error Hamiltonians are combinations of the error Hamiltonian~(\ref{eq:leapfrog:H:err}) for the corresponding leapfrog integrator restricted to $V=V_{\!c}$, the error Hamiltonian~(\ref{eq:BDB:H:err}) of the map \change{[DB]$^2$} restricted to $V=V_{\!s}$, and the mixed term $\{\opn T,V_{\!s}\cls,V_{\!c}\}$. Interestingly, the amplitude of the mixed term is not the same between these four methods: that for map \change{[DKB]$^2$}~(\ref{eq:map:[DKB]^2}) is four times larger than for the maps \change{[DBK]$^2$} and \change{[BKD]$^2$} (\ref{eq:maps:ex:leapfrog}a,\change{b}). Moreover, the maps~(\ref{eq:maps:ex:leapfrog}a,\change{d}) require only one ordinary kick operation $\Exp{h\op{V}_{\!c}}$ per step (either in the middle or at beginning and end, when the accelerations computed in the previous step can be recycled), while the maps~(\ref{eq:maps:ex:leapfrog}b,\change{c}) require two kicks per step. Hence, of the maps~(\ref{eq:maps:ex:leapfrog}) the best computational efficiency to accuracy relation is achieved by the map \change{[DBK]$^2$}~(\ref{eq:map:[DBK]^2}) and the worst by the map \change{[DKB]$^2$}~(\ref{eq:map:[DKB]^2}).

We now verify equation~\eqref{eq:H:err:[KDB]^2} by explicitly monitoring $\tilde{H}_2$, the surrogate Hamiltonian $\tilde{H}$ up to second order (i.e.\ $H$ plus the expressions given in equation~\ref{eq:H:err:[KDB]^2} computed via formul\ae~\ref{eq:VijTT}-\ref{eq:TVijVik}), for an integration of the outer Solar system (the same as in Section~\ref{sec:HB15:4}). We consider three cases: either $\mathcal{S}$ is empty (when the integrator is the ordinary kick-drift-kick leapfrog), $\mathcal{S}$ contains the four Sun-planet pairs, or $\mathcal{S}$ contains all 10 pairs (when the map is identical to \change{[DB]$^2$ = }HB15).

\begin{figure}
	\begin{center}
	\includegraphics[width=0.9\columnwidth]{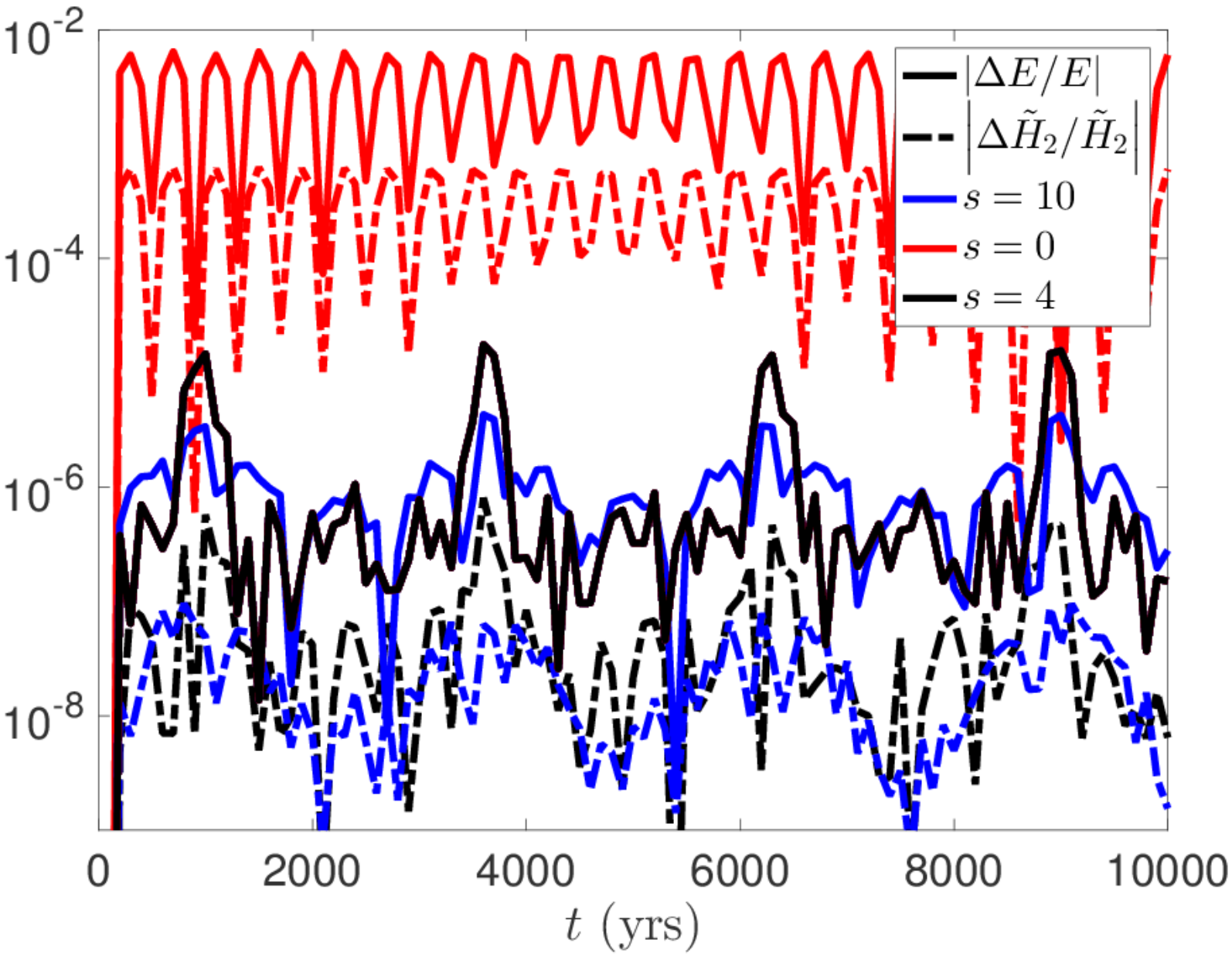}
	\end{center}
	\caption{Verification of equation~\eqref{eq:H:err:[KDB]^2}. $E$ and $\tilde{H}_2$ (the surrogate Hamiltonian up to order $h^2$) are calculated as a function of time for an integration of the outer Solar system using the integrator~\eqref{eq:map:[KDB]^2} for three choices of the set $\mathcal{S}$ (empty: $s=0$, only Sun-planet pair: $s=4$, all pairs: $s=10$). As expected $|\Delta\tilde{H}_2/\tilde{H}_2|$ is smaller than $| \Delta E/E|$ in all cases.
  \label{fig:eq33b}
  }
\end{figure}

Fig.~\ref{fig:eq33b} shows $|\Delta E/E|$ and $|\Delta \tilde{H}_2/\tilde{H}_2|$ as function of time. For all three cases, $|\Delta\tilde{H}_2/\tilde{H}_2|$ is smaller than $|\Delta E/E|$ and $|\Delta\tilde{H}_2/\tilde{H}_2|\propto h^4$ (not shown), confirming equation~\eqref{eq:H:err:[KDB]^2}.

\subsection{A fourth-order hybrid integrator}\label{sec:mix:4}
In order to preserve the fourth-order nature of the integrator, we apply the method of equation~(\ref{eq:CC:4}) and construct the map
\begin{equation} \label{eq:map:[KDBK]^2}
\change{
	\psi_h^{\mathrm{[KDBK]^2}} =\;\;
}
	\Exp{\frac{h}{6}\op{V}_{\!c}}\;
\change{
	\phi_{h/2}^{\dag s}\;
}
	\Exp{\frac{2h}{3}\op{V}_{\!c}}\;
\change{
	\phi_{h/2}^{s\phantom{\dag}}\;
}
	\Exp{\frac{h}{6}\op{V}_{\!c}},
\end{equation}
which combines the maps~(\ref{eq:map:CC:2}) and~(\ref{eq:map:[DB]^2}). In appendix~\ref{app:H:err:mix:4}, we derive the error Hamiltonian of this map to be
\begin{equation} \label{eq:[KDBK]^2:H:err}
	\sub{H}{err}^{\change{\mathrm{[KDBK]^2}}}
	= \frac{h^2}{48}\{\opn T,V_{\!s}\cls,V_{\!s}\}_3
	+ \frac{h^2}{72}\{\opn T,V_{\!c}\cls,V_{\!c}\} + \mathcal{O}(h^4),
\end{equation}
which in the limits of empty set $\mathcal{S}$ or \change{its complement} $\mathcal{S}^c$ obtains the respective previous cases~(\ref{eq:CC:2:H:err}) and (\ref{eq:BDB:H:err}), as expected. Interestingly, the mixed term $\{\opn T,V_{\!s}\cls,V_{\!c}\}$, accounting for three-body interactions \change{$(i,j,k)$ with $(i,j)\in\mathcal{S}$ and $(i,k)\not{\in}\mathcal{S}$}, does not appear. The terms
\begin{equation} \label{eq:Gs:Gc}
	G_s \equiv \{\opn T,V_{\!s}\cls,V_{\!s}\}_3
	\qquad\text{and}\qquad
	G_c \equiv \{\opn T,V_{\!c}\cls,V_{\!c}\}
\end{equation}
can be integrated (see Appendix~\ref{app:G}) to obtain the fourth-order map
(\change{with a parameter $\alpha$ as in equation~\ref{eq:map:[DB]^2_4}})
\begin{eqnarray} 
	\label{eq:map:[KDBK]^2_4}
	\psi_h^{\mathrm{\change{[KDBK]^2_4}}} &=&
	\Exp{\frac{h}{6}\op{V}_{\!c}}\;
	\Exp{-\change{\alpha}\frac{h^3}{96}\op{G}_{\!s}}\;
\change{
	\phi_{h/2}^{\dag s}\;
}
	\Exp{\frac{2h}{3}(\op{V}_{\!c}-\frac{h^2}{48}\op{G}_{\!c})}\;
\change{
	\Exp{(\alpha-1)\frac{h^3}{48}\op{G}_{\!c}}
}
	\nonumber \\ & \cdot &
\change{
	\phi_{h/2}^{s\phantom{\dag}}\;
}
	\Exp{-\change{\alpha}\frac{h^3}{96}\op{G}_{\!s}}\;
	\Exp{\frac{h}{6}\op{V}_{\!c}}.
\end{eqnarray}
\change{This is a generalisation of the fourth-order map~\eqref{eq:map:[DB]^2_4} insofar as it obtains that map when all particle pairs are in set $\mathcal{S}$, and we use `DH16' for both forms. Conversely, when set $\mathcal{S}$ is empty the map~\eqref{eq:map:[KDBK]^2_4}} reduces to the integrator~\change{\eqref{eq:CC:4}}.

We first numerically verify the order of the integrator~\eqref{eq:map:[KDBK]^2_4} \change{with $\alpha=1$}, using the outer Solar system with the four Sun-planet pairs placed in set $\mathcal{S}$ and the six planet-planet pairs in $\mathcal{S}^c$. This grouping is more efficient than the others we tested as we will see below. We integrate for $t=1000$ years and plot in Fig.~\ref{fig:orderd} the absolute energy error at the end of the integration against the step size $h$. The errors are well fit by a $|\Delta E/E| \sim h^4$ curve, as expected, as long as the errors are dominated by truncation (rather than round-off) error. 

\begin{figure}
	\begin{center}
	\includegraphics[width=0.85\columnwidth]{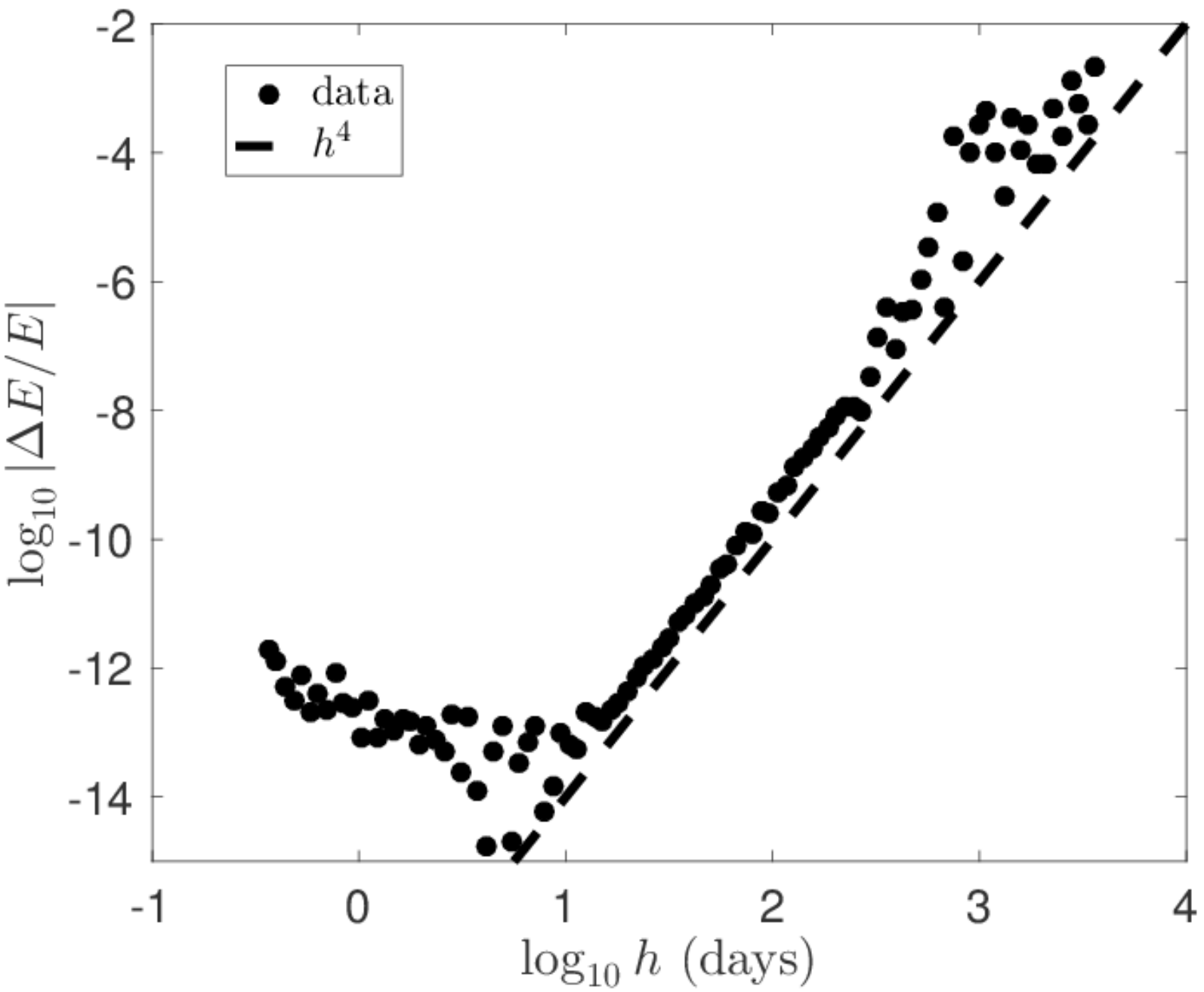}
	\end{center}
	\caption{Verification of the order of map~\eqref{eq:map:[KDBK]^2_4} \change{with $\alpha=1$}. We integrate the outer Solar system with the four Sun-planet pairs in $\mathcal{S}$ for 1000 years. The errors are well fit by $h^4$, except below $h\sim1\,$day, where round-off errors dominate).
  \label{fig:orderd}}
\end{figure}

\subsubsection{Testing conservation of first integrals}
\begin{figure}
	\begin{center}
	\includegraphics[width=0.85\columnwidth]{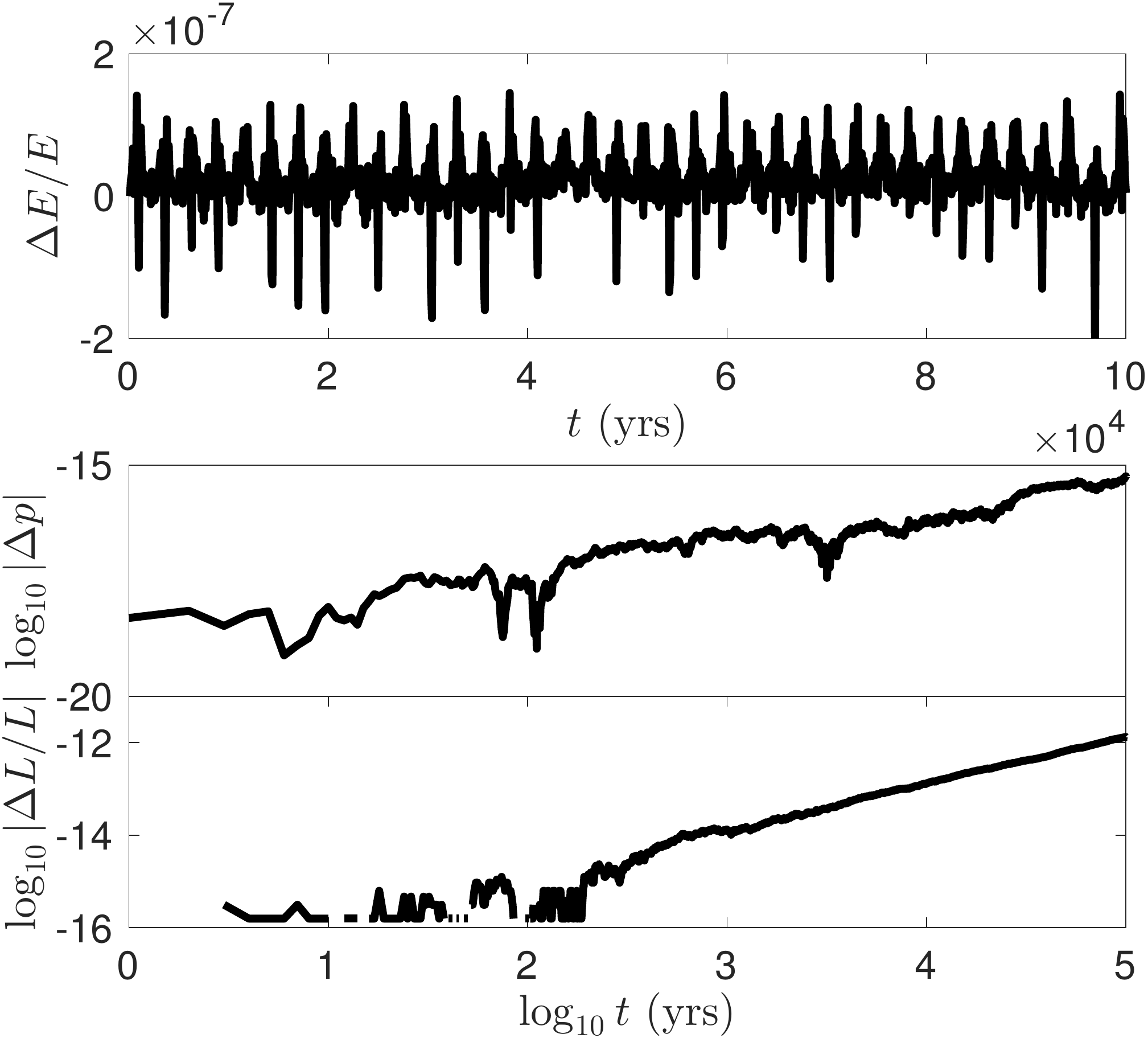}
	\end{center}
	\caption{Conservation of isolating integrals for the integrator~\eqref{eq:map:[KDBK]^2_4} \change{with $\alpha=1$} applied to an integration of the outer Solar system. The energy error is bounded, as expected, while the momentum and angular momentum errors are at the machine precision initially and grow in time due to accumulation of roundoff error.
	\label{fig:outer4vst}
	}
\end{figure}

Next we consider the conservation of isolating integrals. Since the integrator~\eqref{eq:map:[KDBK]^2_4} is a composition of maps that each conserve linear and angular momentum, so does the integrator as a whole. Additionally, the existence of the function $\tilde{H}$ guarantees that the energy error is bounded over exponentially long times \citep{HairerLubichWanner2006}. We test these predictions by integrating the outer Solar system over 100,000 years in steps of $h=1$\,yr. The error in isolating integrals as a function of time is shown in Fig.~\ref{fig:outer4vst}. 
There is no secular drift of the energy error as a function of time as expected. The errors in linear and angular momentum are not exactly zero, but grow like $\propto t^{0.8}$ and $\propto t$, respectively. This is steeper than $t^{1/2}$ expected for accumulation of (unbiased) round-off errors and indicative of bias in the rounding behaviour \citep{Henrici1962}\change{, though our computations use the common IEEE\,754 standard for floating-point arithmetic}. In this test, we used the Kepler solver described by \cite{WisdomHernandez2015}, which shows some bias in tests of a two-body orbit. While such bias can be controlled by careful numerical implementation \citep{ReinTamayo2015}, the value of keeping integration errors near the machine precision, especially for chaotic problems, is questionable \citep{PortegiesZwartBoekholt2014, Hernandez2016}.   

\subsubsection{The effect of set $\mathcal{S}$ on efficiency}
\label{sec:DH16:eff}
As mentioned previously, we are interested in the efficiency of map \eqref{eq:map:[KDBK]^2_4} when the grouping of particle pairs into $\mathcal{S}$ is varied. For three settings of $\mathcal{S}$ (those used in Section~\ref{sec:mix:2} and Fig.~\ref{fig:eq33b} above), we find (by trial and error) the computation time $\sub{t}{cpu}$ required to reach $\langle|\Delta E/E|\rangle\lesssim3\times10^{-8}$. The results are shown in Table \ref{tab:solargroupingsfourth}. 
\begin{table}
	\begin{center}
	\caption{Efficiency of the map~\eqref{eq:map:[KDBK]^2_4} when integrating the outer Solar system for 1000 years for different choices of $\mathcal{S}$. For the first and last choice the method reduces to maps~\eqref{eq:CC:4} and \eqref{eq:map:[DB]^2_4}, respectively.
	\label{tab:solargroupingsfourth}
	}
	\begin{tabular}{llll}
		\hline
		\multicolumn{1}{c}{$\mathcal{S}$}
			& \multicolumn{1}{c}{$\mathcal{S}^c$}
			& \multicolumn{1}{c}{$\langle|\Delta E/E|\rangle$}
			& \multicolumn{1}{c}{$\sub{t}{cpu}$} \\
		\hline
		empty & all & $2.7\times 10^{-8}$ & 0.90\,sec \\
		Sun-planet & planet-planet & $3.0\times 10^{-8}$ & 0.40\,sec \\
		all & empty & $3.1\times 10^{-8}$ & 0.73\,sec \\
		\hline
	\end{tabular}
	\end{center}
\end{table}
The grouping with only the Sun-planet interactions integrated by a Kepler solver is the most efficient, followed by the maps~\eqref{eq:map:[DB]^2_4} (all ten interactions performed via a Kepler solver) and \eqref{eq:CC:4} (no Kepler solver used). One may try to explain this behaviour by studying the error Hamiltonian of the map~\eqref{eq:map:[KDBK]^2_4} up to order $h^4$, but the difficulties in doing so are likely to exceed those we encountered above for the simpler map~\eqref{eq:map:[DB]^2_4}.

\subsubsection{Efficiency comparison with other methods}
\begin{figure}
	\includegraphics[width=0.9\columnwidth]{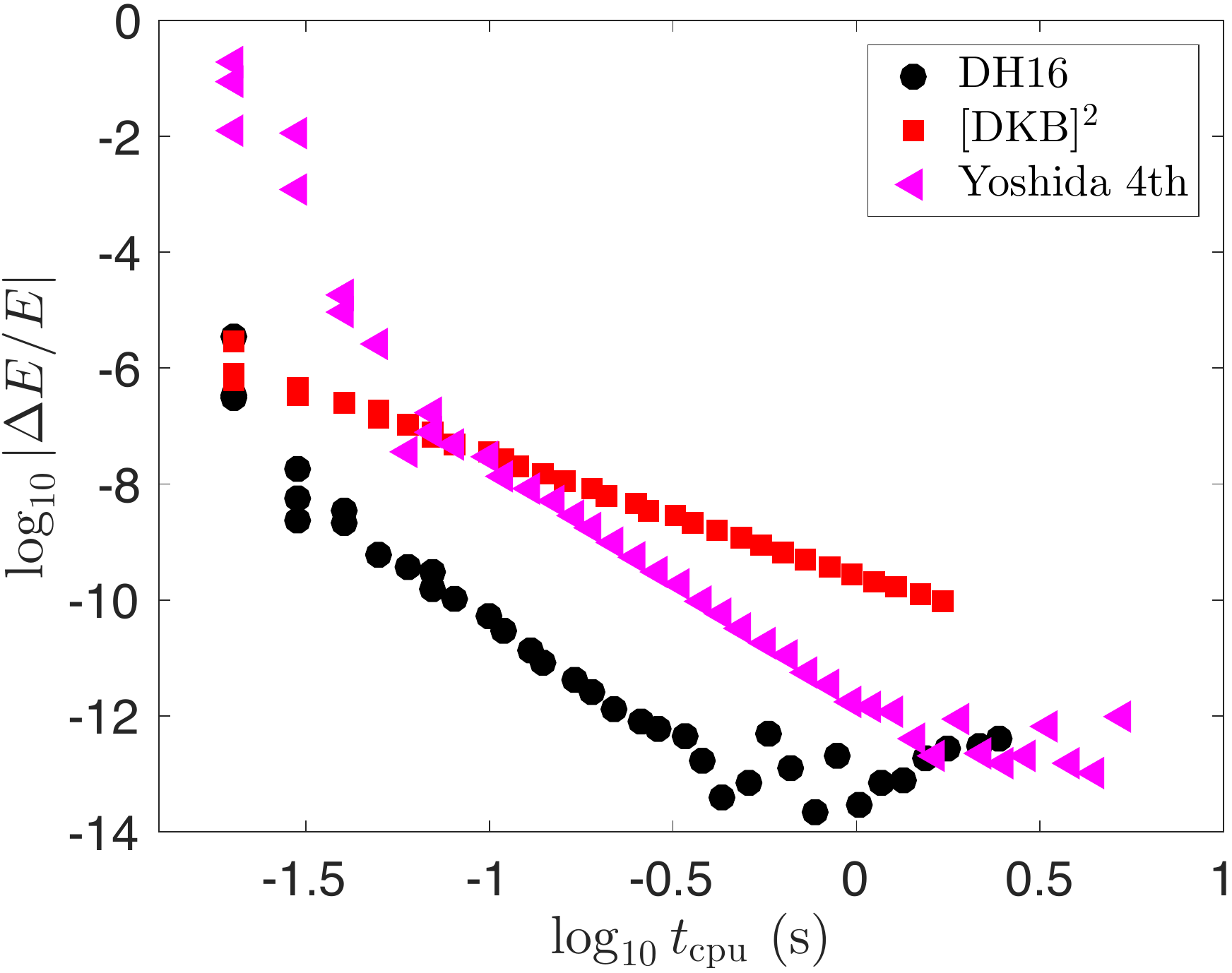}
	\caption{Comparison of efficiency of various integrators (see text) when integrating the outer Solar system for 10,000 years. At high accuracy, our new fourth-order map~\eqref{eq:map:[KDBK]^2_4} (labelled DH16) is the most efficient, but at low accuracy the second-order method \change{[DKB]$^2$} becomes competitive.
\label{fig:methodefficiency}
	}
\end{figure}

Another \change{point} of interest is how map~\eqref{eq:map:[KDBK]^2_4} compares with other integrators in terms of efficiency. We use an integration of the outer Solar system to compare three methods that use Kepler solvers: our new fourth-order map~\change{[KDBK]$^2_4$=\;}DH16 (equation~\ref{eq:map:[KDBK]^2_4}), the second-order method \change{[DKB]$^2$}~(\ref{eq:map:[DKB]^2}), which has the worst efficiency to accuracy relation of the extended leapfrog maps~(\ref{eq:maps:ex:leapfrog}), and a fourth order map obtained by composing three \change{[DKB]$^2$} maps with the recipe of \cite{Yoshida1990}, labelled `Yoshida 4th'. \change{[DKB]$^2$} was compared against other integrators before \citep{HernandezBertschinger2015, Hernandez2016}, and shown to often be the most efficient \change{out of a set of seven published} method\change{s}. \change{DH16, [DKB]$^2$, and ÔYoshida 4thÕ} require a choice for $\mathcal{S}$, and we use the same grouping as in Fig.~\ref{fig:outer4vst}, i.e.\ that for which map~\eqref{eq:map:[KDBK]^2_4} is most efficient. The result is shown in Fig.~\ref{fig:methodefficiency}.      
Map~\eqref{eq:map:[KDBK]^2_4} is most efficient for most of the parameter region shown. In particular, it is always better than the other fourth-order map tested (`Yoshida 4th'). This is not very surprising, since that latter method requires three times as many calls to the Kepler solver.

\section{Discussion and Conclusion}
We have analysed novel symplectic and time reversible integrators for collisional $N$-body problems, where close encounters play an important role in driving the dynamics. These encounters render collisional $N$-body problems much harder than collision-less dynamics and are the main stumbling block for efficient time integration.

Most symplectic integrators which have been applied to collisional dynamics in the past are only second-order accurate and generally handle close encounters inaccurately. A promising approach to overcome this hurdle is the usage of a Kepler solver to deal with close encounters \citep{GoncalvesFerrariEtAl2014}. \cite{HernandezBertschinger2015} have demonstrated how to use this approach to build a symplectic and time-reversible integrator (HB15 \change{or [DB]$^2$ in our nomenclature}). We provide theoretical justification for the success of HB15 and \change{some} related methods: terms of the error Hamiltonian that originate from close two-body encounters are eliminated at all orders. This leaves only close encounters of three or more particles to contribute to the truncation error.

The lowest-order error Hamiltonian of the resulting integration methods can be expressed as the nested Poisson bracket $\{\{T,V\},V\}$ of kinetic and potential energies excluding terms of the form $\{\{T,V_{\!\ij}\},V_{\!\ij}\}$, which account for two-body encounters and are eliminated owing to the Kepler solver ($V_{\!\ij}$ denotes the potential energy arising form the gravitational interaction of particles $i$ and $j$, see equation~\ref{eq:T,V}). Since $T$ is quadratic in the momenta and $V$ a function of the positions only, the term
$\{\{T,V\},V\}$ itself depends only on the particle positions. As a consequence, this terms acts like a potential energy and is integrable. Thus, the associated truncation error can be corrected in a symplectic way and with little extra cost (compared to the solutions of the Kepler problems), resulting in the fourth-order symplectic and time-reversible integrator \change{[DB]$^2_4$} presented in Section~\ref{sec:HB15:4}.

The usage of a Kepler solver may be restricted to a sub-set $\mathcal{S}$ of all pair-wise particle interactions \citep{Hernandez2016}, when the terms $\{\{T,V_{\!\ij}\},V_{\!\ij}\}$ and $\{\{V_{\!\ij},T\},T\}$ from interactions $(i,j)\not{\in}\mathcal{S}$ contribute to the error Hamiltonian. This may be tolerable if such interactions are never close (for example, those between the gas giant planets in the Solar system). However, these terms can also be eliminated in a different way, namely using the method of \cite{Chin1997} which cancels $\{\{V_{\!\ij},T\},T\}$ and integrates $\{\{T,V_{\!\ij}\},V_{\!\ij}\}$ without the need for backward steps (as opposed to the fourth-order symplectic method of \citealt{Yoshida1990}), \change{resulting} in the new symplectic integrator \change{`DH16' of equation~\eqref{eq:map:[KDBK]^2_4}}. \change{This map} is a hybrid between the fourth-order forward method of \cite{Chin1997} and our fourth-order extension \change{[DB]$^2_4$} of HB15\change{, which is its limiting case when all particle pairs are in set $\mathcal{S}$}. 

Various tests and efficiency comparisons of the maps we discuss are presented. As our tests revealed, the novel fourth-order integrators are generally more efficient than previous methods when high accuracy is demanded. However, they still suffer inaccuracies, in particular in some chaotic systems. For a chaotic restricted three-body exchange orbit test \change{and a $N=1024$ cluster simulation} \change{our} fourth-order integrator \change{DH16 with all particle pairs treated with a Kepler solver} performed similarly to \change{the second-order methods} HB15 \change{or (for the cluster simulation only) \textsc{sakura} of \citeauthor{GoncalvesFerrariEtAl2014}, which also used a Kepler solver for each particle pair although in a way that destroys symplecticity and reversibility}. \change{The dynamics of} these \change{systems is} likely dominated by three-body encounters, and the only way to increase the accuracy in such situations appears some form of \change{adaption either of the} time stepping \change{or of the set $\mathcal{S}$ of particle pairs for which a Kepler solver is used. These methods change from one surrogate Hamiltonian to another and may lose symplecticity but retain time reversibility.} We plan to explore \change{these ideas} in the future.

\section*{Acknowledgements}
We thank Edmund Bertschinger for feedback and discussions, Jack Wisdom for suggesting the chaotic exchange orbit used in section~\ref{sec:chaotex}\change{, Guilherme Gon\c{c}alves Ferrari for providing data plotted in Fig~\ref{fig:nbody}, and the anonymous referee for useful comments.

This work used the DiRAC Complexity system, operated by the University of Leicester IT Services, which forms part of the STFC DiRAC HPC Facility (\href{http://www.dirac.ac.uk}{www.dirac.ac.uk}). This equipment is funded by BIS National E-Infrastructure capital grant ST/K000373/1 and  STFC DiRAC Operations grant ST/K0003259/1. DiRAC is part of the National E-Infrastructure}.
\bibliographystyle{mnras}
\bibliography{doc}
\onecolumn
\appendix

\section{Second order error Hamiltonians} \label{app:H:err}
Applying the Campbell-Baker-Haussdorff formula~(\ref{eq:BCH}) twice and trice, we find
\begin{eqnarray}
	\label{eq:BCH:XYX:2}
	\log\big(\Exp{\change{\frac{1}{2}}X}\Exp{Y}\Exp{\change{\frac{1}{2}}X}\big) &=& 
	\change{X}+Y
	- \tfrac{1}{\change{24}}[X,\change{[}X,Y\change{]}]
	+ \tfrac{1}{\change{12}}[Y,\change{[}Y,X\change{]}] + \dots,
	\\
	\label{eq:BCH:XYZYX:2}
	\log\big(\Exp{\change{\frac{1}{2}}X}\Exp{\change{\frac{1}{2}}Y}\Exp{Z}
			 \Exp{\change{\frac{1}{2}}Y}\Exp{\change{\frac{1}{2}}X}\big) &=& 
	\change{X}+\change{Y}+Z
	- \tfrac{1}{\change{24}}[X+Y,\change{[}X+Y,Z\change{]}]
	- \tfrac{1}{\change{24}}[X,\change{[}X,Y\change{]}]
	+ \tfrac{1}{\change{12}}[Z,\change{[}Z,X+Y\change{]}]
	+ \tfrac{\change{1}}{\change{12}}[Y,\change{[}Y,X\change{]}] + \dots.
\end{eqnarray}
Because of equation~(\ref{eq:CommutePoisson}), these relations translate directly to corresponding relations for the surrogate Hamiltonian of a composite map:\begin{eqnarray}
	\label{eq:Herr:ABA:2}
	\Exp{\frac{h}{2}\op{A}}\map{h}{B}\Exp{\frac{h}{2}\op{A}}
	& \quad\text{has}\quad & \tilde{H} = A+B
	- \tfrac{1}{24} h^2 \{\opn B,A\cls,A\} + \tfrac{1}{12} h^2 \{\opn A,B\cls,B\}
	+ \mathcal{O}(h^4)
	\\
	\label{eq:Herr:ABCBA:2}
	\Exp{\frac{h}{2}\op{A}}\Exp{\frac{h}{2}\op{B}}
		\map{h}{C}\Exp{\frac{h}{2}\op{B}}\Exp{\frac{h}{2}\op{A}}
	& \quad\text{has}\quad & \tilde{H} = A+B+C
	- \tfrac{1}{24}h^2\{\opn C,B\cls,B\}
	- \tfrac{1}{24}h^2\{\opn B+C,A\cls,A\}
	+ \tfrac{1}{12}h^2\{\opn B,C\cls,C\}
	+ \tfrac{1}{12}h^2\{\opn A,B+C\cls,B+C\}
	+ \mathcal{O}(h^4).
\end{eqnarray}
\subsection{The method of Hernandez \& Bertschinger}
	\label{app:H:err:HB15}
In order to derive the surrogate Hamiltonians of the schemes~(\ref{eq:map:HB15}), we specify the order in which the maps~(\ref{eq:binary:kick}) are applied in equation~(\ref{eq:map:W}). To this end, we index the $K\equiv N(N-1)/2$ particle pairs
\begin{equation}
	\change{(}i_n,j_n\change{)},\quad n=1\dots K,
\end{equation}
with the implication that pair $(i_1,j_1)$ comes first in the map $\psi_{h}^{W}$ (and last in its adjoint $\psi_{h}^{\dag W}$). If we further define $V_n\equiv V_{i_n\!j_n}$, the map~(\ref{eq:map:[BD]^2}) can be expressed recursively as
\begin{equation} \label{eq:map:[BD]^2:recursive}
	\psi_h^{\mathrm{\change{[BD]^2}}} = \bar{\psi}_h^K,
	\qquad
	\bar{\psi}_h^n =
	\Exp{\frac{h}{2}(\op{V}_n+\op{T})}\;\Exp{-\frac{h}{2}\op{T}}\;\bar{\psi}_h^{n-1}\;
	\Exp{-\frac{h}{2}\op{T}}\;\Exp{\frac{h}{2}(\op{V}_n+\op{T})},
	\qquad
	\bar{\psi}_h^0 = \map{h}{T},
\end{equation}
where we made use of 
\begin{equation} \label{eq:map:Wij:T}
	\Exp{h(\op{V}_{\ij}+\op{T}_i+\op{T}_{\!j})}\,  \Exp{-h(\op{T}_i+\op{T}_{\!j})}
	= \Exp{h(\op{V}_{\ij}+\op{T})}\, \Exp{-h\op{T}_{\phantom{ij}}},
\end{equation}
which follows from equation~(\ref{eq:CommutePoisson}) and $\{T-T_i-T_j,V_{\ij}\}=0$. Applying equation~(\ref{eq:Herr:ABCBA:2}) to the recursion (\ref{eq:map:[BD]^2:recursive}) we find the following recursion for the surrogate Hamiltonian of $\bar{\psi}_h^n$
\begin{eqnarray} \label{eq:map:[BD]^2:H:surr:a}
	\tilde{H}_n &=& \tilde{H}_{n-1} + V_n
	-\frac{h^2}{12}\{\opn T,\tilde{H}_{n-1}\cls,\tilde{H}_{n-1}\}
	+\frac{h^2}{12}\{\opn V_n+T,\tilde{H}_{n-1}-T\cls,\tilde{H}_{n-1}-T\}
	-\frac{h^2}{24}\{\opn \tilde{H}_{n-1},T\cls,T\}
	-\frac{h^2}{24}\{\opn \tilde{H}_{n-1}-T,V_n+T\cls,V_n+T\}
	+\mathcal{O}(h^4)
\end{eqnarray}
and $\tilde{H}_0=T$. For future reference, it proves useful to consider the particular form  $\tilde{H}_n=T + W_n$ where $\partial W_n/\partial\vec{p}_i=\change{\mathcal{O}(h^2)}$. In this case, equation~(\ref{eq:map:[BD]^2:H:surr:a}) reduces to
\begin{eqnarray} \label{eq:map:[BD]^2:H:surr:b}
	\tilde{H}_n &=& T + W_{n-1} + V_n
	+\frac{h^2}{24}\{\opn T,V_n\cls,W_{n-1}\}
	+\mathcal{O}(h^4).
\end{eqnarray}
For the Hernandez \& Bertschinger integrator, we make the ansatz $\tilde{H}_n=T+\bV_n+h^2E_n$ \change{with $E_n$ to be determined and}
\begin{equation}
\change{
	\bV_n \equiv \sum_{k<n} V_k.
}
\end{equation}
\change{E}quation~(\ref{eq:map:[BD]^2:H:surr:b}) \change{then} gives \change{for $W_n=\bV_n+h^2E_n$}
\begin{equation}
	\label{eq:HB15:E2:n}
	E_{n} - E_{n-1} =
	\tfrac{1}{24}\{\opn T,V_n\cls,\bVx\}
	\;+\;\mathcal{O}(h^2).
\end{equation}
\change{Note that $\partial W_n/\partial\vec{p}_i=\mathcal{O}(h^2)$ (as required for equation~\ref{eq:map:[BD]^2:H:surr:b}) follows by induction from $W_0=0$ and the recursion~(\ref{eq:HB15:E2:n}).}

The error Hamiltonian of the complete map then follows as
\begin{equation} \label{eq:BDB:H:err:A}
	\sub{H}{err}^{\mathrm{\change{[BD]^2}}}
	= h^2 E_K + \mathcal{O}(h^4) = \frac{h^2}{24}
		\sum_{n=1}^K \{\opn T,V_n\cls,\bVx\} + \mathcal{O}(h^4)
	= \frac{h^2}{24} \sum_{n=1}^{K}
		\sum_{k=1}^{n-1} \{\opn T,V_{i_nj_n}\cls,V_{i_kj_k}\} + \mathcal{O}(h^4).
\end{equation}
The double sum in this last form includes each pair of pairs $\{\opn T,V_{\ij}\cls,V_{lk}\}$ exactly once, except for those with $(i,j)=(l,k)$ which are not contained at all. The term $\{\opn T,V\cls,V\}$, which occurs in the error Hamiltonian of the leapfrog, contains each pair of pairs $\{\opn T,V_{\ij}\cls,V_{lk}\}$ twice, except for $(i,j)=(l,k)$ which are  contained once. Thus, the form~(\ref{eq:BDB:H:err:A}) is identical to
\begin{equation} \label{eq:H:err:BDB}
	\sub{H}{err}^{\mathrm{\change{[BD]^2}}}
	= \frac{h^2}{48}\big(\{\opn T,V\cls,V\}-\{\opn T,V\cls,V\}_2\big)
	+ \mathcal{O}(h^4)
	= \frac{h^2}{48}\{\opn T,V\cls,V\}_3 + \mathcal{O}(h^4).
\end{equation}
Writing \change{(with $H_n=V_n+T$)}
\change{
\begin{eqnarray}
	\phi_h^{\dag} =
	\Exp{ h\op{T}}\,
	\psi_{h}^{\dag W_1}\;\cdots\;\psi_{h}^{\dag W_K} &=&
	\Exp{ h\op{T}}\;
	\Exp{-h\op{T}}\Exp{h\op{H}_1}\;\cdots\;
	\Exp{-h\op{T}}\Exp{h\op{H}_K}\,
	\nonumber \\ &=&
	\label{eq:map:[DB]:recursive}
	\phantom{\Exp{ h\op{T}}\,\Exp{-h\op{T}}\,}
	\Exp{ h\op{H}_1}\;\cdots\;
	\Exp{-h\op{T}}\,\Exp{h\op{H}_K}\,
	\Exp{-h\op{T}}\,
	\Exp{ h\op{T}}
	=\;
	\psi_{h}^{W_1}\;\cdots\;\psi_{h}^{W_K}\,
	\Exp{ h\op{T}},
\end{eqnarray}
}
we see that \change{$\phi_h^{\dag}$ and $\phi_h=\psi_{h}^{W_K}\cdots\psi_{h}^{W_1}\,\Exp{h\op{T}}$ are identical} except for \change{a reversal of the order of} binary kicks.

\begin{figure}
	\begin{minipage}[b]{85mm}
		\includegraphics[width=70mm]{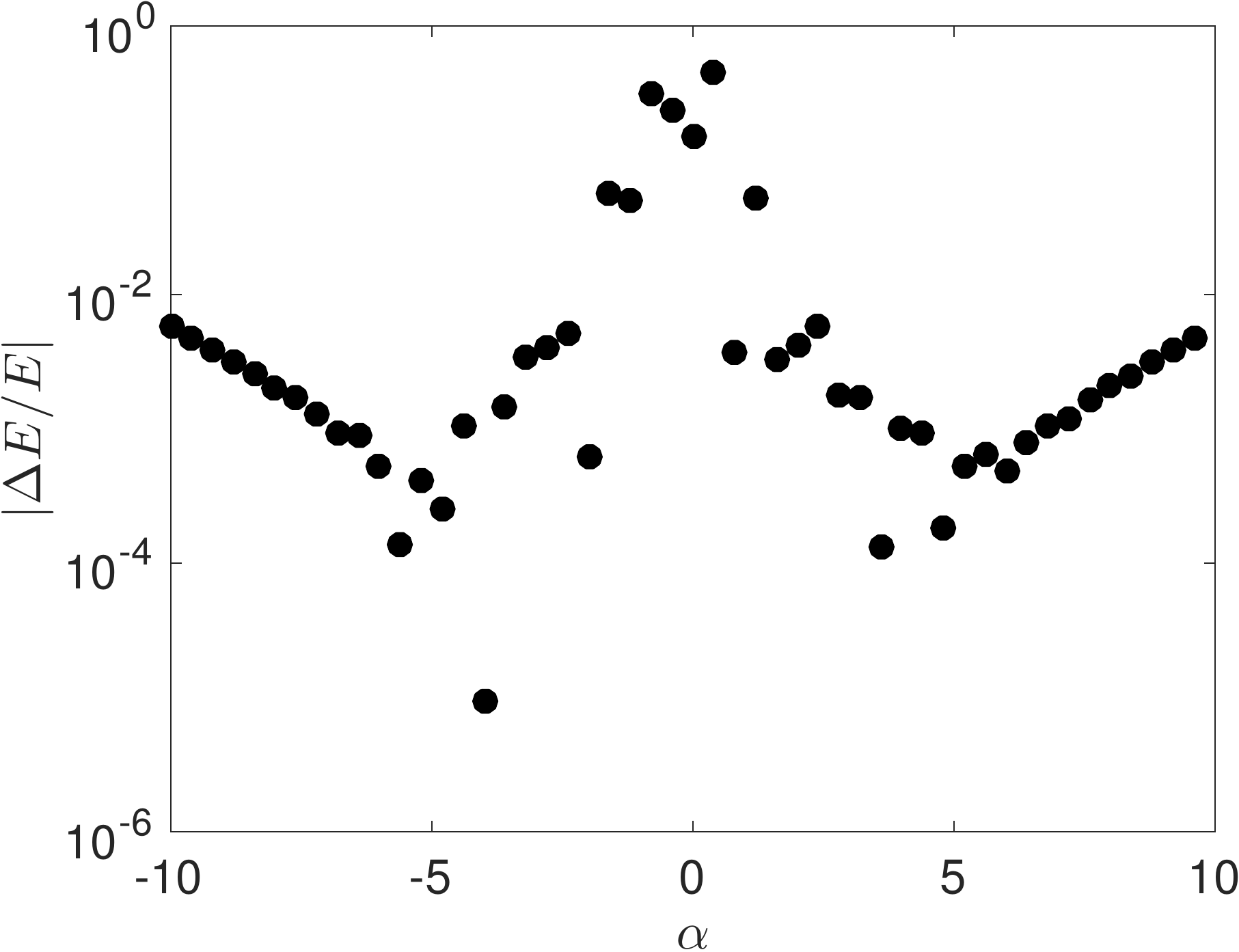}
	\end{minipage}
	\hfill
	\begin{minipage}[b]{85mm}
	\caption{Energy error as a function of $\alpha$ for map \eqref{eq:map:psi:alpha} when integrating an elliptic Kepler orbit with eccentricity $e=0.9$ over one in 100 steps. Only $\alpha=0$ (not shown) gives an error at the level of machine precision.
	\label{fig:alpha29}
	}
	\end{minipage}
\end{figure}

The error Hamiltonian of the map~(\ref{eq:map:psi:alpha}) can be obtained analogously
\change{to that of map~\eqref{eq:map:[BD]^2}} as
\begin{eqnarray}
	\label{eq:err:map:psi:alpha}
	\sub{H}{err}^\alpha &=&
	\frac{\alpha(\alpha-1)h^2}{6}\{\opn V,T\cls,T\}
	+\frac{\alpha h^2}{4}\{\opn T,V\cls,V\}_2
	+\frac{(1+12\alpha)h^2}{48}\{\opn T,V\cls,V\}_3
	\;+\;\mathcal{O}(h^4).
\end{eqnarray}
We demonstrate the error properties of the map~(\ref{eq:map:psi:alpha}) by integrating the equal-mass two-body problem with elliptic orbit of eccentricity $e=0.9$ over one period in 100 equal time steps.  For $\alpha=0$, the magnitude of the energy error approaches the computational round-off error, while for all other values the energy error becomes substantial as shown in Fig.~\ref{fig:alpha29}. 
Some $\alpha$ are better than others, but the smallest error is approximately $10^{-5}$.  Thus it is essential to let $\alpha=0$.

\subsection{Error Hamiltonian for the extended Leapfrog}
\label{app:H:err:mix:2}
\change{The map [BDK]$^2$, i.e.\ $\phi_{h/2}^{s\phantom{\dag}}\Exp{h\op{V}_{\!c}}\phi_{h/2}^{\dag s}$, differs from~\eqref{eq:map:[DBK]^2} only in the order of binary kicks, which as we will see has no effect on the second-order error terms.}
The maps \change{[BDK]$^2$} and~(\ref{eq:map:[BKD]^2}) have the same recursive form as the map~(\ref{eq:map:[BD]^2}), but start from $\bar{\psi}_h^0=\Exp{\frac{h}{2}\op{T}}\,\Exp{h\op{V}_{\!c}}\,\Exp{\frac{h}{2}\op{T}}$ and $\bar{\psi}_h^0=\Exp{\frac{h}{2}\op{V}_{\!c}}\,\Exp{h\op{T}}\,\Exp{\frac{h}{2}\op{V}_{\!c}}$, respectively. Consequently, the recursion for the respective surrogate Hamiltonian is identical to equation~(\ref{eq:map:[BD]^2:H:surr:a}), except that
\begin{eqnarray}
	\tilde{H}_0 &=& T+V_{\!c} - \frac{h^2}{24}\{\opn V_{\!c},T\cls,T\}
		+ \frac{h^2}{12}\{\opn T,V_{\!c}\cls,V_{\!c}\} + \mathcal{O}(h^4)
	\qquad\text{for $\psi^{\mathrm{\change{[BDK]^2}}}$ and}
	\\
	\tilde{H}_0 &=&
	T+V_{\!c} + \frac{h^2}{12}\{\opn V_{\!c},T\cls,T\}
		- \frac{h^2}{24}\{\opn T,V_{\!c}\cls,V_{\!c}\} + \mathcal{O}(h^4)
	\qquad\text{for $\psi^{\mathrm{\change{[BKD]^2}}}$.}
\end{eqnarray}
With the ansatz $\tilde{H}_n=T+V_{\!c} + \bV_n + h^2 E_n + \mathcal{O}(h^4)$, we obtain from equation~(\ref{eq:map:[BD]^2:H:surr:b})
\begin{equation}
	E_{n} - E_{n-1} =
	 \tfrac{1}{24}\{\opn T,V_n\cls,\bVx\}
	+\tfrac{1}{24}\{\opn T,V_n\cls,V_{\!c}\}
	\;+\;\mathcal{O}(h^2)
\end{equation}
and therefore
\begin{eqnarray}
	\tilde{H}^{\mathrm{\change{[BDK]^2}}} &=& 
	\tilde{H}_K = T+V
	- \frac{h^2}{24}\{\opn V_{\!c},T\cls,T\}
	+ \frac{h^2}{24}\{\opn T,V_{\!s}\cls,V_{\!c}\}
	+ \frac{h^2}{12}\{\opn T,V_{\!c}\cls,V_{\!c}\}
	+ \frac{h^2}{48}\{\opn T,V_{\!s}\cls,V_{\!s}\}_3
	+ \mathcal{O}(h^4),
	\\
	\tilde{H}^{\mathrm{\change{[BKD]^2}}} &=&
	\tilde{H}_K = T+V
	+ \frac{h^2}{12}\{\opn V_{\!c},T\cls,T\}
	+ \frac{h^2}{24}\{\opn T,V_{\!s}\cls,V_{\!c}\}
	- \frac{h^2}{24}\{\opn T,V_{\!c}\cls,V_{\!c}\}
	+ \frac{h^2}{48}\{\opn T,V_{\!s}\cls,V_{\!s}\}_3
	+ \mathcal{O}(h^4).
\end{eqnarray}

The error Hamiltonian~(\ref{eq:H:err:[KDB]^2}) for the integrator~(\ref{eq:map:[KDB]^2}) follows directly from that of the integrator~(\ref{eq:map:[BD]^2}) derived in equation~(\ref{eq:H:err:BDB}) above and equation~(\ref{eq:Herr:ABA:2}) to account for the maps $\Exp{\frac{h}{2}\op{V}_{\change{\!c}}}$ at the beginning and end.

The map~(\ref{eq:map:[DKB]^2}) can be written
\begin{equation}
	\psi_h^{\mathrm{\change{[DKB]^2}}} =
	\Exp{\frac{h}{2}\op{T}}\;
	\Exp{\frac{h}{2}\op{V}_{\!c}}\;
	\Exp{-\frac{h}{2}\op{T}}\;
	\psi_h^{\mathrm{\change{[DB]^2}}}\;
	\Exp{-\frac{h}{2}\op{T}}\;
	\Exp{\frac{h}{2}\op{V}_{\!c}}\;
	\Exp{\frac{h}{2}\op{T}},
\end{equation}
when threefold application of equation~(\ref{eq:Herr:ABA:2}) starting from $\tilde{H}^{\mathrm{\change{[DB]^2}}}$ gives the error Hamiltonian reported in equation~(\ref{eq:H:err:[DKB]^2}).

\subsection{Error Hamiltonian for the integrator of section~\ref{sec:mix:4}}\label{app:H:err:mix:4}
\change{The map~(\ref{eq:map:[KDBK]^2}) differs from
\begin{equation} \label{eq:map:[KBDK]^2}
	\psi_h^{\mathrm{[KBDK]^2}} =\;\;
	\Exp{\frac{h}{6}\op{V}_{\!c}}\;
\change{
	\phi_{h/2}^{s\phantom{\dag}}\;
}
	\Exp{\frac{2h}{3}\op{V}_{\!c}}\;
\change{
	\phi_{h/2}^{\dag s}\;
}
	\Exp{\frac{h}{6}\op{V}_{\!c}}
\change{
	\;\;=
	\Exp{\frac{h}{6}\op{V}_{\!c}}\;
	\psi_{h/2}^{W_s}\;
	\Exp{\frac{h}{2}\op{T}}\;
	\Exp{\frac{2h}{3}\op{V}_{\!c}}\;
	\Exp{\frac{h}{2}\op{T}}\;
	\psi_{h/2}^{\dag W_s}\;
	\Exp{\frac{h}{6}\op{V}_{\!c}}
},
\end{equation}
only in the order of binary kicks, because of equation~\eqref{eq:map:[DB]:recursive}.} We rewrite the map~(\ref{eq:map:[KBDK]^2}) recursively as
\begin{eqnarray} \label{eq:[KDBK]^2:recursive}
	\Exp{\frac{h}{6}\op{V}_{\!c}}\,\bar{\psi}_h^K\,\Exp{\frac{h}{6}\op{V}_{\!c}},
	\qquad
	\bar{\psi}_h^n =
	\Exp{\frac{h}{2}(\op{V}_n+\op{T})}\,\Exp{-\frac{h}{2}\op{T}}\,\bar{\psi}_h^{n-1}\,
	\Exp{-\frac{h}{2}\op{T}}\,\Exp{\frac{h}{2}(\op{V}_n+\op{T})},
	\qquad
	\bar{\psi}_h^0 = \Exp{\frac{h}{2}\op{T}}\,\Exp{\frac{2h}{3}\op{V}_{\!c}}
		\,\Exp{\frac{h}{2}\op{T}}
\end{eqnarray}
where again $V_n\equiv V_{i_nj_n}$ for $(i_n,j_n)\in\mathcal{S}$. The recursion relation for $\bar{\psi}_h^n$ defined in~(\ref{eq:[KDBK]^2:recursive}) only differs from that of equation~(\ref{eq:map:[BD]^2:recursive}) by the starting point. Consequently, the recursion for its surrogate Hamiltonian is identical to equation~(\ref{eq:map:[BD]^2:H:surr:a}), but with 
\begin{equation}
	 \tilde{H}_0 = T + \frac{2}{3}V_{\!c}
	 -\frac{h^2}{36}\{\opn V_{\!c},T\cls,T\} + \frac{h^2}{27}\{\opn T,V_{\!c}\cls,V_{\!c}\}
	 +\mathcal{O}(h^4).
\end{equation}
With the ansatz $\tilde{H}_n=T+\tfrac{2}{3}V_{\!c}+\bV_n+h^2E_n+\mathcal{O}(h^4)$, we obtain from equation~(\ref{eq:map:[BD]^2:H:surr:b})
\begin{equation}
	E_{n}-E_{n-1} = 
	 \tfrac{1}{36}\{\opn T,V_n\cls,V_{\!c}\}
	+\tfrac{1}{24}\{\opn T,V_n\cls,\bVx\}
	\;+\;\mathcal{O}(h^2)
\end{equation}
and therefore
\begin{equation}
	\tilde{H}_K = T + V_{\!s} + \frac{2}{3}V_{\!c}
	- \frac{h^2}{36}\{\opn V_{\!c},T\cls,T\}
	+ \frac{h^2}{27}\{\opn T,V_{\!c}\cls,V_{\!c}\}
	+ \frac{h^2}{36}\{\opn T,V_{\!s}\cls,V_{\!c}\}
	+ \frac{h^2}{24}\sum_{n=1}^K\{\opn T,V_n\cls,\bVx\}
	\;+\;\mathcal{O}(h^4)	
\end{equation}
Finally, the surrogate and error Hamiltonian of the complete map~(\ref{eq:map:[KBDK]^2}) follows from one last application of equation~(\ref{eq:Herr:ABA:2}) to account for the maps $\Exp{\frac{h}{6}\op{V}_{\!c}}$ at begin and end
\begin{equation} \label{eq:hyberror}
	\sub{H}{err}^{\mathrm{\change{[KBDK]^2}}}
	= \frac{h^2}{48}\{\opn T,V_{\!s}\cls,V_{\!s}\}_3
	+ \frac{h^2}{72}\{\opn T,V_{\!c}\cls,V_{\!c}\} + \mathcal{O}(h^4),
\end{equation}
in particular the mixed term $\{\opn T,V_{\!s}\cls,V_{\!c}\}$ does not appear. 
\change{Since at second order $\sub{H}{err}^{\mathrm{[KBDK]^2}}$ does not depend on the order of binary kicks, $\sub{H}{err}^{\mathrm{[KDBK]^2}}=\sub{H}{err}^{\mathrm{[KBDK]^2}}+\mathcal{O}(h^4)$.}

\section{Fourth-order Error Hamiltonians}
The Campbell-Baker-Haussdorff formula~(\ref{eq:BCH}) up to order five reads
\begin{eqnarray}\label{eq:BCH:5}
	\log\big(\Exp{X}\Exp{Y}\big) &=& 
	X+Y + \tfrac{1}{2}[XY] + \tfrac{1}{12}\big([X^2Y] + [Y^2X]\big)
	-\tfrac{1}{24}[YX^2Y]
	\nonumber \\ && \phantom{X+Y}
	-\tfrac{1}{720}\big([X^4Y]+[Y^4X]\big)
	+\tfrac{1}{360}\big([XY^3X]+[YX^3Y]\big)
	+\tfrac{1}{120}\big([XYXYX]+[YXYXY]\big)
	\;\dots,
\end{eqnarray}
where we have used a compact bracket notation, e.g.\ $[XY^3X]=[X,Y,Y,Y,X]=[X,[Y,[Y,[Y,X]]]]$. With this, we can extend equations~(\ref{eq:BCH:XYX:2}) and 
(\ref{eq:Herr:ABA:2}) to fourth order:
\begin{eqnarray}
	\label{eq:BCH:XYX:4}
	\log\big(\Exp{\change{\frac{1}{2}}X}\Exp{Y}\Exp{\change{\frac{1}{2}}X}\big) &=& 
	\change{X}+Y - \tfrac{1}{\change{24}}[X^2Y] + \tfrac{1}{\change{12}}[Y^2X]
	+\tfrac{7}{\change{5760}}[X^4Y]
	-\tfrac{1}{\change{720}}[Y^4X]
	+\tfrac{1}{\change{360}}[YX^3Y]
	+\tfrac{1}{\change{360}}[XY^3X]
	+\tfrac{1}{\change{120}}[YXYXY]
	-\tfrac{1}{\change{480}}[XYXYX] \dots,
\end{eqnarray}
\begin{eqnarray}
	\label{eq:Herr:ABA:4}
	\tilde{H} &=& \change{A}+B
	+ \tfrac{1}{\change{12}} h^2 \{AB^2\}
	- \tfrac{1}{\change{24}} h^2 \{BA^2\}
	+ \tfrac{7}{\change{5760}} h^4 \{BA^4\}
	- \tfrac{1}{\change{720}} h^4 \{AB^4\}
	+ \tfrac{1}{\change{360}} h^4 \{AB^3A\}
	+ \tfrac{1}{\change{360}} h^4 \{BA^3B\}
	- \tfrac{1}{\change{480}} h^4 \{ABABA\}
	+ \tfrac{1}{\change{120}} h^4 \{BABAB\}
	+\mathcal{O}(h^6)
\end{eqnarray}
(e.g.\ \citealt{Yoshida1990}, equation 3.2; \citealt{HairerLubichWanner2006}, equation 4.15), using the compact notation also for Poisson brackets. With this relation we can compute the error Hamiltonian of any self-adjoint composite symplectic map to fourth order.

\subsection{The fourth-order error terms}
\label{app:Err:4}
Note that the nested Poisson brackets $\{TV^4\}$ and $\{TV^3T\}$ vanish (regardless of their coefficients). The remaining fourth-order error terms can be split into contributions from two-, three-, and four-body encounters: $\{VT^3V\} = \{VT^3V\}_2 + \{VT^3V\}_3$, $\{TVTVT\} = \{TVTVT\}_2 + \{TVTVT\}_3$, and $\{VTVTV\} = \{VTVTV\}_2 + \{VTVTV\}_3 + \{VTVTV\}_4$. Using the notation $V_n\equiv V_{i_n\!j_n}$ from Appendix~\ref{app:H:err:HB15}, we have
\begin{subequations}
\begin{eqnarray}
	\{VT^3V\}_2
	&=& \sum_n \{V_nT^3V_n\},
	\\
	\{VT^3V\}_3
	&=& \sum_{n}\sum_{k\neq n}\{V_nT^3V_k\},
	\\
	\{TVTVT\}_2
	&=& \sum_n \{TV_nTV_nT\},
	\\
	\{TVTVT\}_3
	&=& \sum_{n}\sum_{k\neq n} \{TV_nTV_kT\},
	\\
	\{VTVTV\}_2
	&=& \sum_{n} \{V_nTV_nTV_n\},
	\\ \label{eq:vtvtv3}
	\{VTVTV\}_3
	&=& \sum_{n}\sum_{k\neq n} \{V_nTV_nTV_k\} + 2\{V_nTV_kTV_k\},
	\\
	\{VTVTV\}_4
	&=& \sum_{n}\sum_{k\neq n}\sum_{l\neq n,k} \{V_nTV_kTV_l\},
\end{eqnarray}
\end{subequations}
where we have used $\{V_nTV_k\}=\{V_kTV_n\}$.
These error terms can be constructed from the following elementary terms, ordered by the number of particles contributing (all indices are distinct).
\begin{subequations}
	\label{eq:Err4:2}
\begin{eqnarray}
	\{V_{\ij}T^4\} &=&
		-3\frac{Gm_im_j}{r_{\ij}^9}
		\left[35(\vec{v}_{\ij}\cdot\vec{x}_{\ij})^4
		 	 -30r_{\ij}^2v_{\ij}^2(\vec{v}_{\ij}\cdot\vec{x}_{\ij})^2
		  	 +3 r_{\ij}^4v_{\ij}^4\right],
	\\
	\{V_{\ij}T^3V_{\ij}\} &=&
		-9\frac{G^2m_im_j(m_i+m_j)}{r_{\ij}^8}
		\left[3(\vec{v}_{\ij}\cdot\vec{x}_{\ij})^2
			 -r_{\ij}^2v_{\ij}^2\right],
	\\
	\{TV_{\ij}TV_{\ij}T\} &=& \phantom{-}
		4 \frac{G^2m_im_j(m_i+m_j)}{r_{\ij}^8}
		\left[6(\vec{v}_{\ij}\cdot\vec{x}_{\ij})^2 - r_{\ij}^2v^2_{\ij}\right],
	\\
	\{V_{\ij}TV_{\ij}TV_{\ij}\} &=&
		-4\frac{G^3m_im_j(m_i+m_j)^2}{r_{\ij}^7};
\end{eqnarray}
\end{subequations}
\begin{subequations}
	\label{eq:Err4:3}
\begin{eqnarray}
	\label{eq:VijTTTVik}
	\{V_{\ij}T^3V_{ik}\} &=&
		-9\frac{G^2m_im_jm_k}{r_{\ij}^7r_{ik}^3}
		\left[
			(\vec{x}_{\ij}\cdot\vec{x}_{ik})
			\big[5(\vec{v}_{\ij}\cdot\vec{x}_{\ij})^2-r_{\ij}^2v_{\ij}^2\big]
			-2r_{\ij}^2(\vec{x}_{\ij}\cdot\vec{v}_{\ij})
						(\vec{x}_{ik}\cdot\vec{v}_{\ij})
		\right],
	\\
	\label{eq:TVijTVikT}
	\{TV_{\ij}TV_{ik}T\} &= &
	6\frac{G^2m_im_jm_k}{r_{\ij}^7r_{ik}^3}
		\left[5(\vec{x}_{\ij}\cdot\vec{v}_{\ij})^2
			   (\vec{x}_{\ij}\cdot\vec{x}_{ik})
			 -{2}r_{\ij}^2(\vec{x}_{\ij}\cdot\vec{v}_{\ij})
			   (\vec{x}_{ik}\cdot\vec{v}_{\ij})
			 -r_{\ij}^2v_{\ij}^2(\vec{x}_{\ij}\cdot\vec{x}_{ik}\})
		\right],
	\nonumber \\ &&+
	2\frac{G^2m_im_jm_k}{r_{\ij}^5r_{ik}^5}
		\left[9(\vec{x}_{\ij}\cdot\vec{v}_{\ij})
			   (\vec{x}_{ik}\cdot\vec{v}_{ik})(\vec{x}_{\ij}\cdot\vec{x}_{ik})
			 -3r_{\ij}^2{(\vec{x}_{ik}\cdot\vec{v}_{\ij})}
			 	(\vec{x}_{ik}\cdot\vec{v}_{ik})
			 -3r_{ik}^2(\vec{x}_{\ij}\cdot\vec{v}_{\ij})
			 	(\vec{x}_{\ij}\cdot\vec{v}_{ik})
			 + r_{\ij}^2r_{ik}^2(\vec{v}_{\ij}\cdot\vec{v}_{ik})
		\right],
	\\
	\label{eq:VijTVijTVik}
	\{V_{\ij}TV_{\ij}TV_{ik}\} &=&
		-4\frac{G^3m_im_jm_k(m_i+m_j)}{r_{\ij}^6r_{ik}^3}
			(\vec{x}_{\ij}\cdot\vec{x}_{ik}),
	\\
	\label{eq:VijTVikTVij}
	\{V_{\ij}TV_{ik}TV_{\ij}\} &=&
		-2\frac{G^3m_im_jm_k(m_i+m_j)}{r_{\ij}^6r_{ik}^3}
			(\vec{x}_{\ij}\cdot\vec{x}_{ik})
		-\frac{G^3m_im_j^2m_k}{r_{\ij}^6r_{ik}^5}
		\left[3(\vec{x}_{\ij}\cdot\vec{x}_{ik})^2-r_{\ij}^2r_{ik}^2\right],
	\\
	\label{eq:VijTVikTVik}
	\{V_{\ij}TV_{ik}TV_{ik}\} &=&
		-2\frac{G^3m_im_jm_k(m_i+m_k)}{r_{\ij}^3r_{ik}^6}
			(\vec{x}_{\ij}\cdot\vec{x}_{ik})
		-\frac{G^3m_im_jm_k^2}{r_{\ij}^5r_{ik}^6}
		\left[3(\vec{x}_{\ij}\cdot\vec{x}_{ik})^2-r_{\ij}^2r_{ik}^2\right];
\end{eqnarray}
\end{subequations}
\begin{subequations}
	\label{eq:Err4:4}
\begin{eqnarray}
	\{V_{\ij}TV_{ik}TV_{\!jk}\} &=& \phantom{-}
		\frac{G^3m_im_jm_k^2}{r_{\ij}^5r_{ik}^3r_{\!jk}^3}
		\left[3(\vec{x}_{\ij}\cdot\vec{x}_{ik})(\vec{x}_{\ij}\cdot\vec{x}_{\!jk})
			-r_{\ij}^2(\vec{x}_{ik}\cdot\vec{x}_{\!jk})\right]
		-\frac{G^3m_im_j^2m_k}{r_{\ij}^3r_{ik}^5r_{\!jk}^3}
		\left[3(\vec{x}_{\ij}\cdot\vec{x}_{ik})(\vec{x}_{ik}\cdot\vec{x}_{\!jk})
			-r_{ik}^2(\vec{x}_{\ij}\cdot\vec{x}_{\!jk})\right],
	\\
	\{V_{\ij}TV_{ik}TV_{il}\} &=&
		-\frac{G^3m_im_jm_km_l}{r_{\ij}^5r_{ik}^5r_{il}^3}
		\left[3([r_{ik}^2\vec{x}_{\ij}+r_{\ij}^2\vec{x}_{ik}]\cdot\vec{x}_{il})
			(\vec{x}_{\ij}\cdot\vec{x}_{ik})
			-r_{\ij}^2r_{ik}^2
			([\vec{x}_{ik}+\vec{x}_{\ij}]\cdot\vec{x}_{il})
		\right],
	\\
	\{V_{\ij}TV_{ik}TV_{jl}\} &=& \phantom{-}
		\frac{G^3m_im_jm_km_l}{r_{\ij}^5r_{ik}^3r_{jl}^3}
		\left[3(\vec{x}_{\ij}\cdot\vec{x}_{ik})(\vec{x}_{\ij}\cdot\vec{x}_{jl})
			-r_{\ij}^2(\vec{x}_{ik}\cdot\vec{x}_{jl})
		\right].
\end{eqnarray}
\end{subequations}
The three-body-encounter terms originating from Poisson brackets with just two $V$ components, (\ref{eq:VijTTTVik}) and (\ref{eq:TVijTVikT}), depend on the particle masses just through the product, while the distance of the first particle pair in each Poisson bracket tends to be more important. In case of the Poisson brackets with three $V$ components, for the three-body-encounter terms
(\ref{eq:VijTVijTVik}), (\ref{eq:VijTVikTVij}), and (\ref{eq:VijTVikTVik}) the masses and distance of the particle pair that appears twice are more important.

\subsection{The Leapfrog integrator}
For the \change{kick-drift-kick} and \change{drift-kick-drift} leapfrog integrators, we obtain immediately from equation~(\ref{eq:Herr:ABA:4})
\begin{eqnarray}
	\sub{H}{err}^{\mathrm{\change{[KD]^2}}} &=&
	- \frac{h^2}{24} \{TV^2\} + \frac{h^2}{12} \{VT^2\}
	- \frac{h^4}{720} \{VT^4\} + \frac{h^4}{120} \{TVTVT\}
	+ \frac{h^4}{360} \{VT^3V\} - \frac{h^4}{480} \{VTVTV\}
	+\mathcal{O}(h^6),
	\\
	\sub{H}{err}^{\mathrm{\change{[DK]^2}}} &=& \phantom{-}
	  \frac{h^2}{12} \{TV^2\} - \frac{h^2}{24} \{VT^2\}
	+ \frac{7h^4}{5760} \{VT^4\} - \frac{h^4}{480} \{TVTVT\}
	+ \frac{h^4}{360} \{VT^3V\} + \frac{h^4}{120} \{VTVTV\}
	+\mathcal{O}(h^6).
\end{eqnarray}

\subsection{The fourth-order extension\change{s} of the method of Hernandez \& Bertschinger}
\label{sec:fourtherror}
The error Hamiltonian of the method~(\ref{eq:map:[DB]^2_4}) presented in section~\ref{sec:HB15:4} differs from that of the original method of \citeauthor{HernandezBertschinger2015} only by the additional terms arising from the correction maps. \change{As in appendix~\ref{app:H:err:HB15} before, we consider the first pair to be  the innermost in the recursive formulation of the integrator. Thus, we actually analyse the map [KBDK]$^2_4$ rather than [KDBK]$^2_4$.} The fourth-order version of equation~(\ref{eq:map:[BD]^2:H:surr:b}) can be calculated from equation~(\ref{eq:Herr:ABA:4}) as
\begin{eqnarray}
	\label{eq:Herr:HB15:5}
	\tilde{H}_n &=& T+W_{n-1}+V_n + \frac{h^2}{24}\{TV_nW_{n-1}\}
	-\frac{h^4}{17280}\{V_nT^3W_{n-1}\}
	-\frac{h^4}{2160}\{W_{n-1}T^3V_n\}
	\nonumber \\ &-&
	 \frac{h^4}{1440}\{TV_nTW_{n-1}T\}
	+\frac{h^4}{1920}\{V_nTV_nTW_{n-1}\}
	+\frac{h^4}{5760}\{V_nTW_{n-1}TV_n\}
	+\frac{h^4}{360}\{V_nTW_{n-1}TW_{n-1}\}
	-\frac{h^4}{480}\{W_{n-1}TW_{n-1}TV_n\}
	+\mathcal{O}(h^6).
\end{eqnarray}
Analogously to the analysis in Appendix~\ref{app:H:err:HB15}, we write the surrogate Hamiltonian of the intermediate maps $\bar{\psi}^n$ as
\begin{equation}
	\label{eq:Hsurr:HB15}
	\tilde{H}_n = T + \bV_n + h^2 E_{2,n} + h^4 E_{4,n} + \mathcal{O}(h^6).
\end{equation}
The second-order term $E_{2,n}$ has contributions from the second-order error~(\ref{eq:HB15:E2:n}) and the correction map $\Exp{(\alpha-1)\frac{h^3}{48}\op{G}}$, i.e.\
\begin{equation}
	\label{eq:Hsurr:HB15:E2:4}
	E_{2,n}
	= \frac{1}{24}\sum_{k<n}\sum_{l<k}\{TV_kV_l\}
	+ \frac{\alpha-1}{24}\sum_{k}\sum_{l<k}\{TV_kV_l\},
\end{equation}
which satisfies $\partial E_{2,n}/\partial\vec{p}_i=0$, and hence commutes with any potential. 
This contributes a second-order error
\begin{equation}
	\label{eq:map:[DB]^2_4:E2}
	h^2 E_{2,K}
	= \frac{\alpha h^2}{24}\sum_n\sum_{k<n}\{TV_nV_k\}
	= \frac{\alpha h^2}{48} \{TVV\}_3.
\end{equation}
The increment of the fourth-order error terms resulting from $W_n=E_{2,n}$ in equation~(\ref{eq:Herr:HB15:5}) is 
\begin{equation}
	\label{eq:map:[DB]^2_4:E4:3}
	(E_{4,n}-E_{4,n-1})_{[3]} = \frac{1}{24}\{TV_nE_{2,n-1}\}
	= \frac{1}{576}\sum_{k<n}\sum_{l<k}\{V_kTV_lTV_n\}
	+ \frac{\alpha-1}{576}\sum_{k}\sum_{l<k}\{V_kTV_lTV_n\}.
\end{equation}
The resulting contribution to the fourth-order error Hamiltonian follows as
\begin{eqnarray}
	\label{eq:map:[DB]^2_4:E4:3:a}
	h^4 E_{4,K[3]} &=&
	  \frac{h^4}{576}\sum_n\sum_{k<n}\sum_{l<k}\{V_kTV_lTV_n\}
	+ \frac{(\alpha-1)h^4}{576}\sum_n\sum_{k}\sum_{l<k}\{V_kTV_lTV_n\}
	\nonumber \\ & &
	\nonumber \\ &=&
	  \frac{h^4}{576}\sum_{n}\sum_{k<n}\sum_{l<k}\{V_kTV_lTV_n\}
	+ \frac{(\alpha-1)h^4}{576}\sum_{n}\sum_{k\neq n} \{V_nTV_kTV_k\}
	+ \frac{(\alpha-1)h^4}{1152}\{VTVTV\}_4,
\end{eqnarray}
where we have re-arranged the sums (as well as re-labeled the indices and exploited $\{V_kTV_n\}=\{V_nTV_k\}$) in order to separate contributions from three- and four-body encounters (we address the overlap between the first two terms later). The correction terms at the beginning and end contribute second- and fourth-order errors {(using equation~(\ref{eq:Herr:ABA:4}) with $A=\alpha h^3\sum_n\sum_{k<n}\{TV_nV_k\}/48$ and $B=T+V$)}
\begin{eqnarray}
	\label{eq:HB15:E2+E4:end}
	&&
	- \frac{\alpha h^2}{24} \sum_{k}\sum_{l<k} \{TV_kV_l\}
	+ \frac{\alpha h^4}{{288}} \sum_{k}\sum_{l<k} \{TV_kV_l(T+V)(T+V)\}
	\nonumber \\ & &
	\nonumber \\ &=&
	- \frac{\alpha h^2}{48} \{TVV\}_3
	+ \frac{\alpha h^4}{{288}} \sum_{n}\sum_{k<n} \{TV_nTV_kT\}
	- \frac{\alpha h^4}{{288}} \sum_{n}\sum_{k\neq n} \{V_nTV_kTV_k\}
	- \frac{\alpha h^4}{{576}} \{VTVTV\}_4
\end{eqnarray}
and the combined contributions to the error Hamiltonian from the terms~(\ref{eq:map:[DB]^2_4:E2}), (\ref{eq:map:[DB]^2_4:E4:3:a}), and (\ref{eq:HB15:E2+E4:end}) is
\begin{eqnarray}
	\label{eq:HB15:E4:3+end}
	  \frac{h^4}{576}\sum_{n}\sum_{k<n}\sum_{l<k}\{V_kTV_lTV_n\}
	- \frac{{(\alpha+1)}h^4}{576} \sum_{n}\sum_{k\neq n} \{V_nTV_kTV_k\}
	- \frac{{(\alpha+1)}h^4}{1152} \{VTVTV\}_4
	+ \frac{\alpha h^4}{{288}} \sum_{n}\sum_{k<n} \{TV_nTV_kT\},
\end{eqnarray}
in particular, the second-order error vanishes (by construction).

The contributions to $E_{4,n}-E_{4,n-1}$ from the fifth-order terms in equation~(\ref{eq:Herr:HB15:5}) are 
\begin{eqnarray}
	\label{eq:HB15:2:E4:5}
	(E_{4,n}-E_{4,n-1})_{[5]} &=& 
	-\tfrac{1}{17280}\{V_nT^3\bVx\}
	-\tfrac{1}{2160}\{\bVx T^3V_n\}
	-\tfrac{1}{1440}\{TV_nT\bVx T\}
	\nonumber \\[0.5ex] &&
	+\tfrac{1}{1920}\{V_nTV_nT\bVx \}
	+\tfrac{1}{5760}\{V_nT\bVx TV_n\}
	+\tfrac{1}{360}\{V_nT\bVx T\bVx \}
	-\tfrac{1}{480}\{\bVx T\bVx TV_n\},
\end{eqnarray}
which contributes
\begin{eqnarray}
	\label{eq:HB15:2:E4:5:a}
	h^4E_{4,K} &=&
	-\frac{h^4}{17280}\sum_{n}\sum_{k<n}\{V_nT^3V_k\}
	-\frac{h^4}{2160} \sum_{n}\sum_{k<n}\{V_kT^3V_n\}
	{-}\frac{h^4}{1440} \sum_{n}\sum_{k<n}\{TV_nTV_kT\}
	\nonumber \\ &&
	+\frac{h^4}{1920}\sum_{n}\sum_{k<n}\{V_nTV_nTV_k\}
	+\frac{h^4}{5760}\sum_{n}\sum_{k<n}\{V_nTV_kTV_n\}
	+\frac{h^4}{360}\sum_{n}\sum_{k<n}\sum_{l<n}\{V_nTV_kTV_l\}
	-\frac{h^4}{480}\sum_{n}\sum_{k<n}\sum_{l<n}\{V_lTV_kTV_n\}
\end{eqnarray}
to the error Hamiltonian. The total error Hamiltonian is the sum of~(\ref{eq:HB15:E4:3+end}) and~(\ref{eq:HB15:2:E4:5:a}). 

\begin{figure}
	\begin{minipage}[b]{85mm}
		\includegraphics[width=80mm]{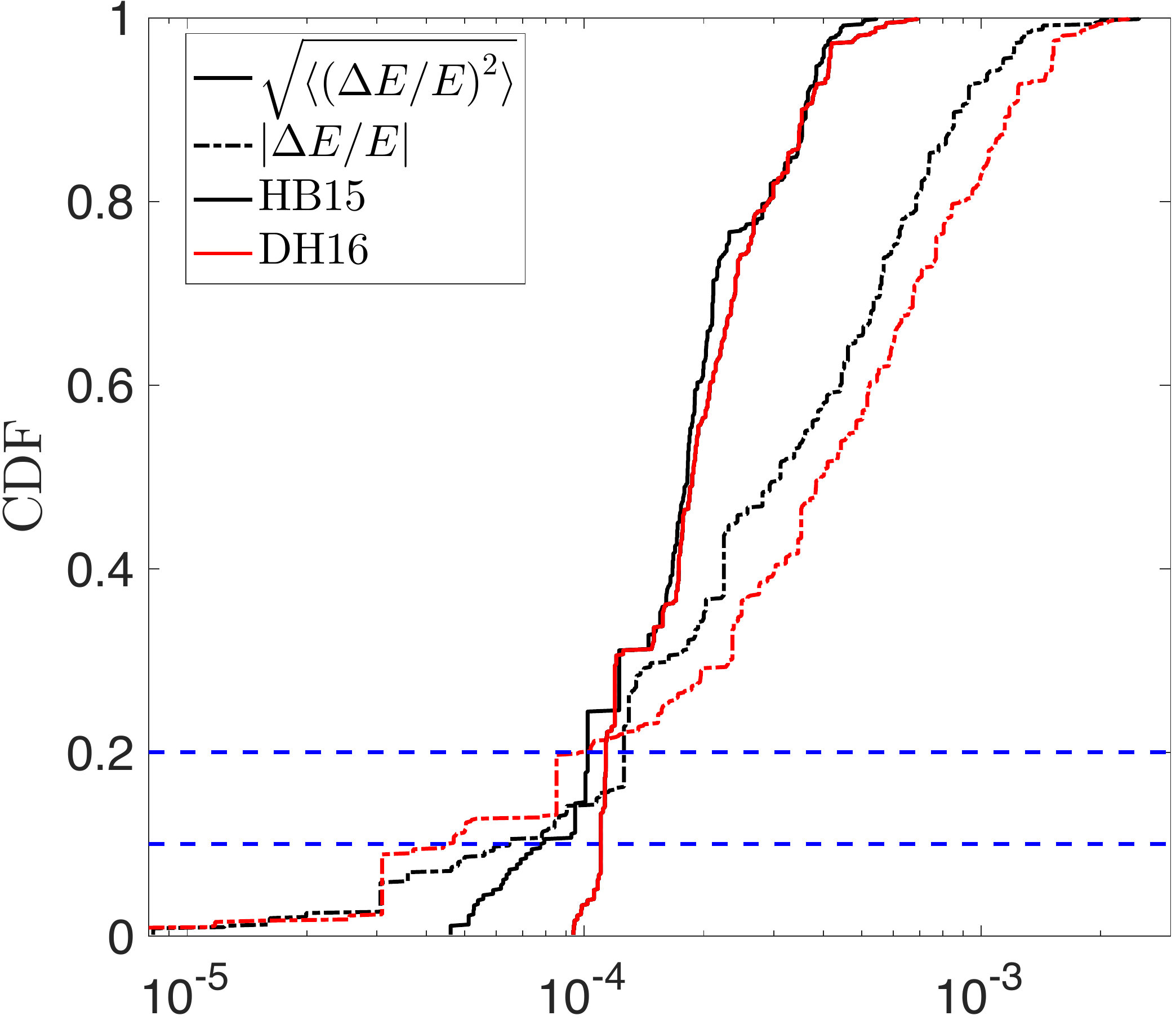}
	\end{minipage}
	\hfill
	\begin{minipage}[b]{85mm}
		\caption{Error distribution for HB15 \change{(=[DB]$^2$)} and its fourth-order accurate extension \change{DH16 (=[DB]$^2_4$~(\ref{eq:map:[DB]^2_4}) with $\alpha=0$)}, for all 720 pair orderings of the hierarchical quadruple problem described in the text. \change{The distribution of rms energy errors is much narrower than that of absolute errors, indicating that the tail of low energy errors is mostly due to chance agreements of the final with the initial total energy. A Kolmogornov-Smirnov analysis suggests a significant difference between the CDFs from the two integrators.}
		\label{fig:orderings}
		}
	\end{minipage}
\end{figure}
In principle, equations~\eqref{eq:HB15:E4:3+end} and \eqref{eq:HB15:2:E4:5:a} in conjunction with equations~\eqref{eq:Err4:3} and \eqref{eq:Err4:4} express the dependence of the error Hamiltonian on the order in which the binary kicks are applied. However, it appears beyond human reasoning to obtain much useful insight from these equations. Therefore, we now explore the effect of the order of binary kicks by numerical experiments.

To this end, a simple problem with widely separated $V_{\ij}$ seems useful. We choose an equal-mass, co-planar, co-rotating, aligned, symmetric, hierarchical quadruple system. The outer equal-mass binary has $e=0.5$ and is initially at apo-centre. Each of its components is in turn an equal mass tighter binary with 100 times smaller semi-major axis and $e=0.9$, one starting from peri-centre, the other from apo-centre (to break degeneracies in the pair potentials). All three binaries are co-planar, co-aligned (the eccentricity vectors point in the same direction), and rotate anti-clockwise; the period ratio between inner and outer binary is $\sqrt{100^3/2}\approx707$. We integrate this system using map~\eqref{eq:map:[DB]^2_4} with $\alpha=0$ for half the period of the outer binary using steps equal to $0.14$ times the inner binary period.  The energy error is bounded in time if $h$ is small enough.

We perform a separate integration for each of the $6!=720$ possible orders of the particle pairs and measure the \change{accumulated} energy error at the end of each integration \change{and the rms energy error over the course of the integration}. Fig.~\ref{fig:orderings} shows the resulting cumulative distribution function\change{s} (CDF\change{s}) together with the equivalent result\change{s} for the integrator HB15 \change{(=[DB]$^2$)}. \change{The horizontal dashed blue lines separate the 10th and 20th percentiles of the CDFs. The rms errors show less variation and thus have smaller tails in their CDFs.} According to our analysis, the ordering does not affect the second-order error terms of HB15, but only its fourth-order errors. For a range of $h$, the magnitude of the errors of the two methods is similar because for this particular problem HB15 behaves similar to a fourth-order method, \change{most likely because three-body encounters (which are solely responsible for the second-order error of HB15) contribute negligibly to the overall error which instead is dominated by four-body encounters (which contribute only at fourth and higher orders).}

\change{We expect the underlying error distributions of the solid lines to be different.  A two-sample KS test supports this expectation and rejects the null hypothesis that the underlying distributions are the same at the 0.1\% level.  The same statements hold for the dashed lines.  The dashed lines show the spread in errors is larger for DH16, but this relative spread disappears in the solid lines.  We ask whether there is a pattern to the orderings corresponding to the low error tail of the dashed red curve: we did not find such a pattern.  We investigated whether the two tight pairs are in a preferential location in the orderings in the best 10\% of the CDF: are the tight pairs usually adjacent, separated, or at the beginning.  The answer to all these questions is no.}

\section{Implementation details}
\subsection{The force gradient terms} \label{app:G}
The map $\Exp{-h^3\op{G}}$ requires a second loop over all particles pairs. After the ordinary accelerations due to $V$,
\begin{equation} \label{eq:acc}
	\vec{a}_{i} = - \frac{1}{m_i} \pdiff{V}{\vec{x_i}} = - \sum_{j\neq i} \frac{Gm_j}{r_{\ij}^3}\vec{x}_{\ij},
\end{equation}
are computed in a first loop, the accelerations due to 
\begin{equation}
	\label{eq:app:TVV}
	G = \{T,V,V\} = \sum_i \frac{1}{m_i} \pdiff{V}{\vec{x_i}}\cdot\pdiff{V}{\vec{x_i}}
	\change{
	\;= \sum_i m_i^{}\, \vec{a}_i^2
	  = - \sum_{i<j} \frac{Gm_im_{\!j}}{r_{\ij}^3}\vec{x}_{\ij}\cdot\vec{a}_{\ij}
	\qquad\qquad\text{with}\qquad
	\vec{a}_{\ij}\equiv\vec{a}_{i}-\vec{a}_{\!j}
	}
\end{equation}
can be computed in a second loop as
\begin{equation} \label{eq:acc:grad}
	\vec{g}_i = - \frac{1}{m_i} \pdiff{G}{\vec{x_i}} = 
	2 \sum_{j\neq i}
	\frac{Gm_j}{r_{\ij}^5}
	\left[\vec{a}_{\ij}\,r_{\ij}^2
		-3 \vec{x}_{\ij}\,(\vec{a}_{\ij}\cdot\vec{x}_{\ij})
	\right].
\end{equation}
The \change{accelerations required for the} map \change{$\Exp{-h^3\op{G}_s}$ and} generated by the term
\begin{equation}
	\label{eq:app:TVV3}
	G_s \;\change{=}\; \{T,V,V\}_3 
	\change{
	\;=\;
	- \sum_{i<j} \frac{Gm_im_{\!j}}{r_{\ij}^3}\vec{x}_{\ij}\cdot\tilde{\vec{a}}_{\ij}
	\qquad\qquad\text{with}\qquad
	\tilde{\vec{a}}_{\ij} \equiv \vec{a}_{\ij} + \frac{G(m_i+m_j)}{r_{\ij}^3}
		\vec{x}_{\ij}
	}
\end{equation}
\change{are} calculated in a similar way \change{as}
\begin{eqnarray}
	\tilde{\vec{g}}_i = -\frac{1}{m_i}\pdiff{G_s}{\vec{x}_i}
	&=& 2 \sum_{j\neq i}
	\frac{Gm_j}{r_{\ij}^5}
	\left[\tilde{\vec{a}}_{\ij}\,r_{\ij}^2
		-3 \vec{x}_{\ij}\,(\tilde{\vec{a}}_{\ij}\cdot\vec{x}_{\ij})
	\right].
\end{eqnarray}
\change{Note that $\tilde{\vec{a}}_{\ij}$ is the difference between the accelerations of particles $i$ and $j$ and due to all other particles, while $\vec{a}_{\ij}$ includes their mutual attraction.}
The accelerations generated by $G_{\!s}$ and $G_{\!c}$ defined in equation~(\ref{eq:Gs:Gc}) are computed analogously, except that only pair-wise interactions contained in, respectively, sets $\mathcal{S}$ and $\mathcal{S}^c$ are considered.

\begin{figure}
	\begin{minipage}[t]{85mm}
		\begin{center}
		\includegraphics[width=40mm]{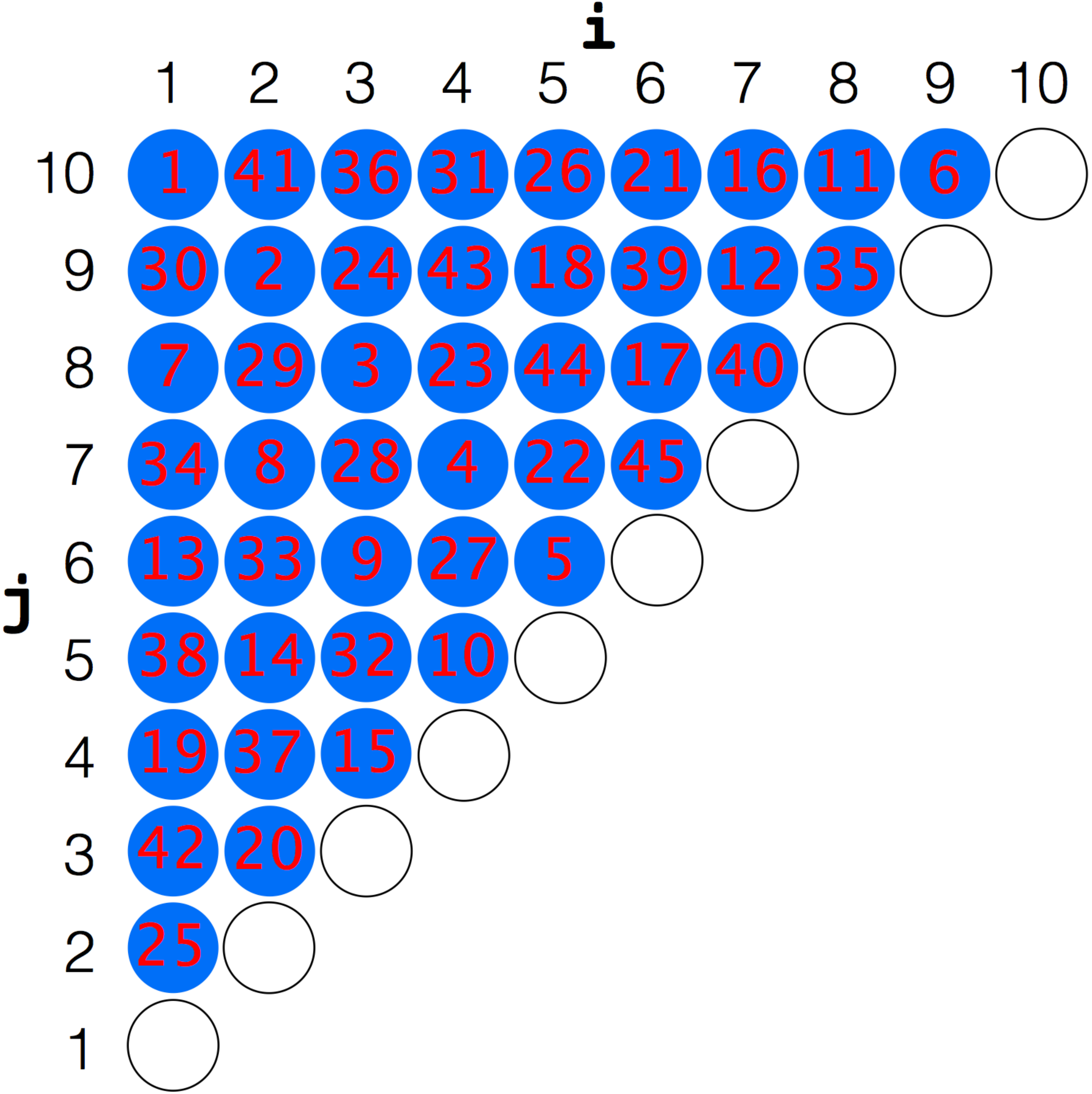}
		\includegraphics[width=40mm]{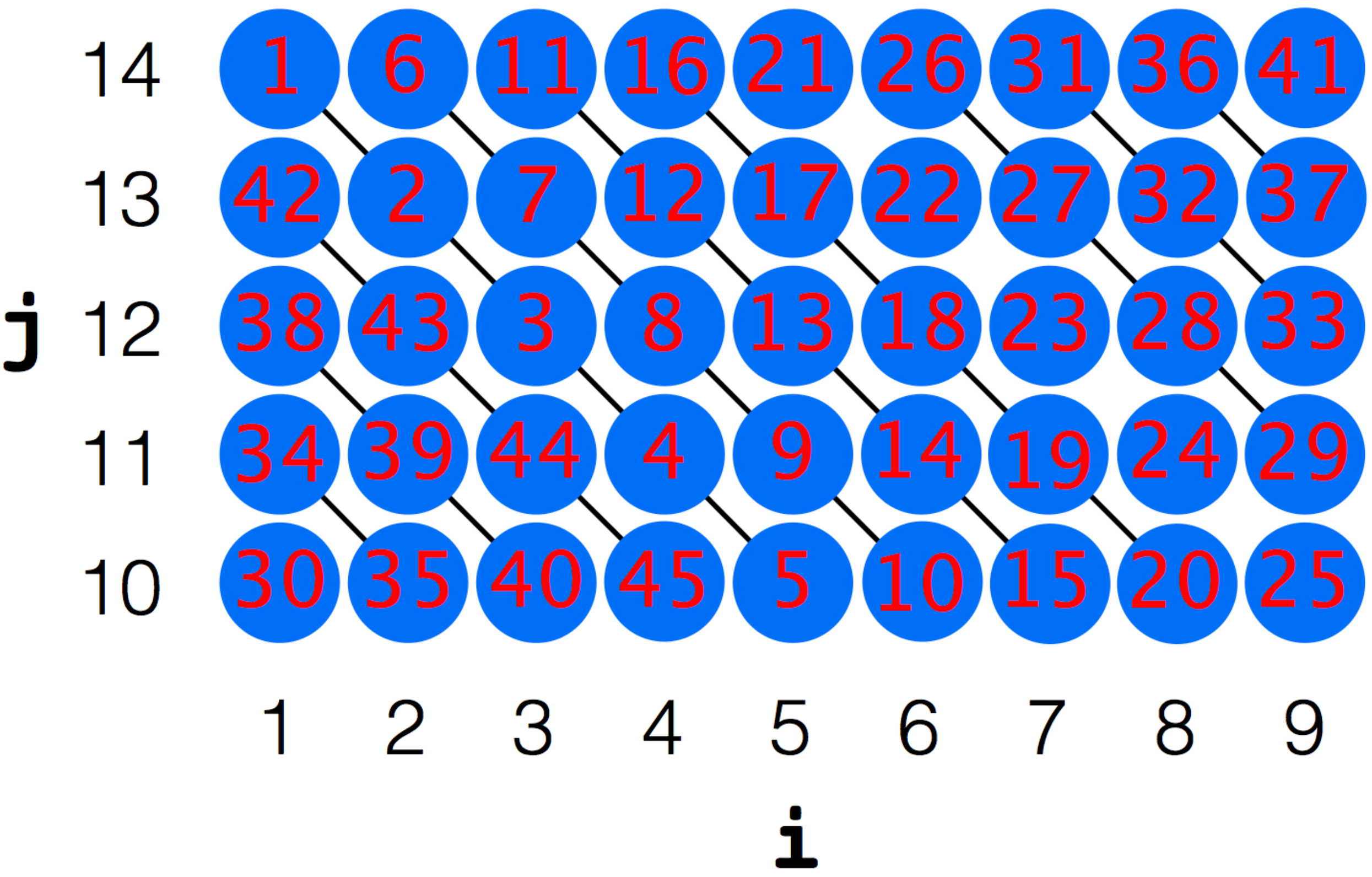}
		\end{center}
		\caption{\change{
Vectorisation of binary kicks must avoid mutually dependent interactions within the same vector. \textbf{Left}: the $N(N-1)/2$ interactions (blue discs) between $N$ (=10 in this example) particles are vectorised in the order indicated in red (obtained by the round-robin method, see text) or its reverse for the adjoint map, and requires vector size $\sub{n}{vec}\le\lfloor N/2\rfloor$. \textbf{Right}: the $N\times M$ interactions between two distinct particle sets are most easily vectorised using a diagonal periodic pattern with vector size $\sub{n}{vec}\le\min(N,M)$.
		\label{fig:vectorise}
		}}
	\end{minipage}
	\hfill
	\begin{minipage}[t]{85mm}
		\begin{center}
		\includegraphics[width=40mm]{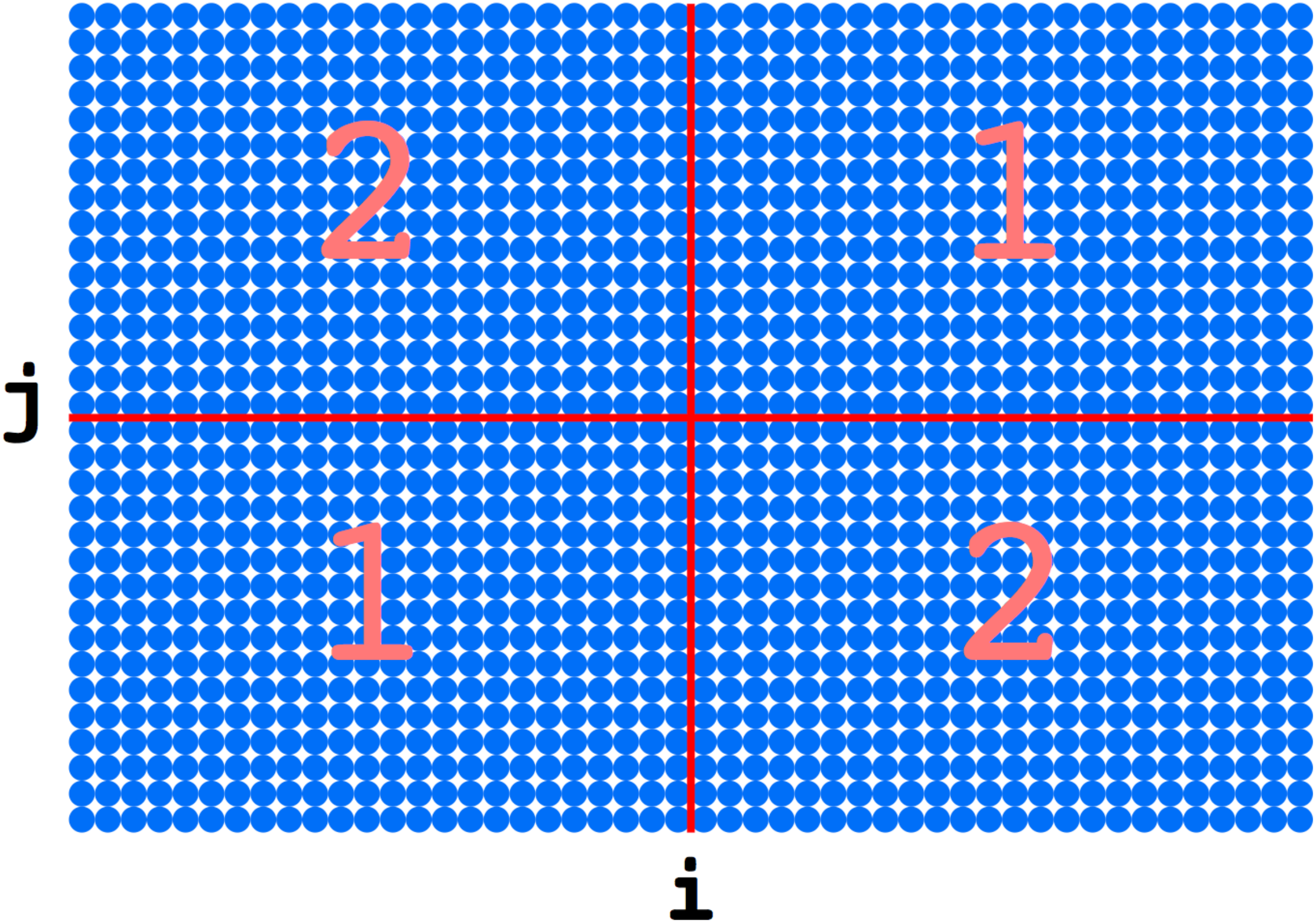}
		\hfill
		\includegraphics[width=40mm]{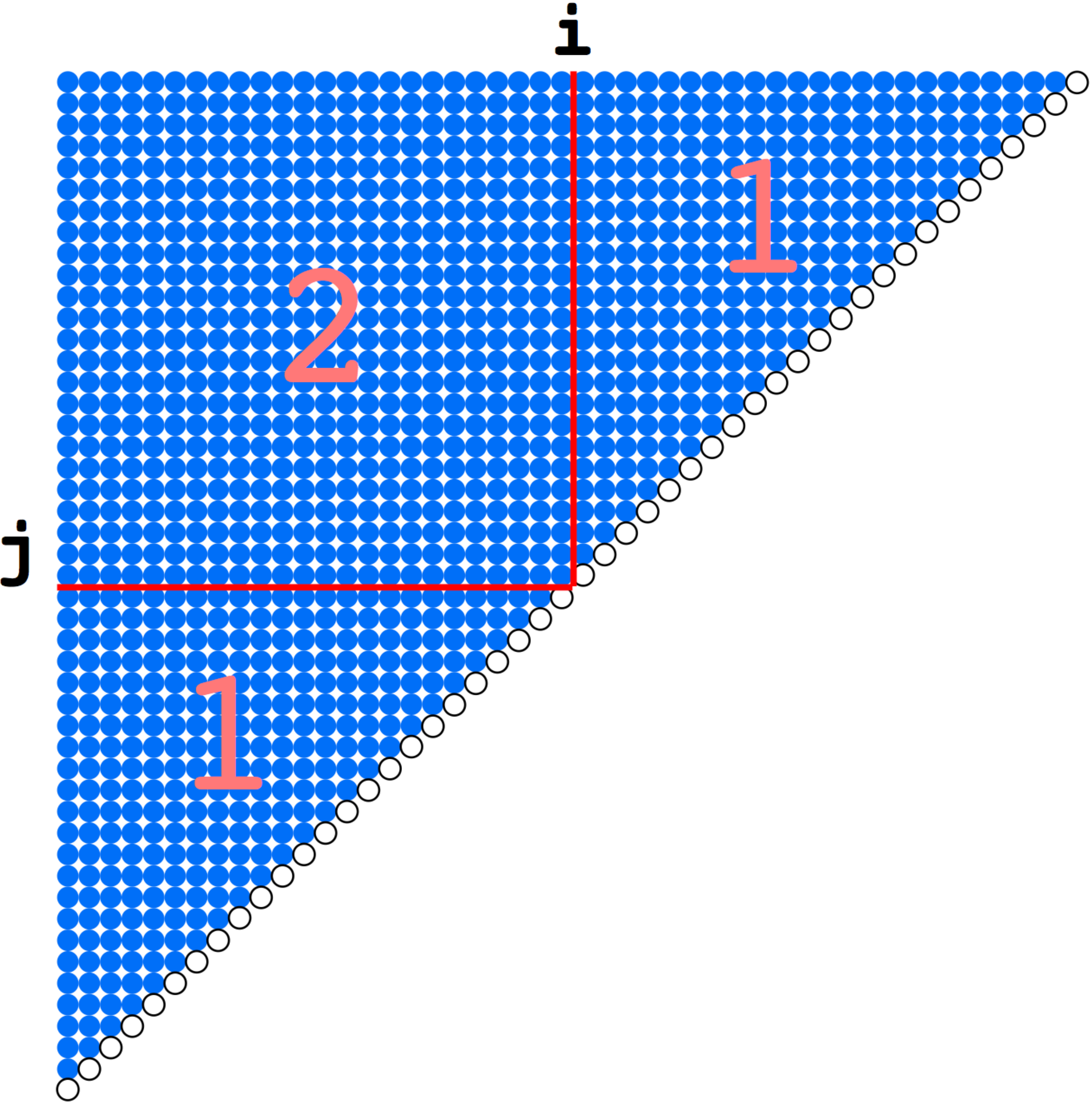}
		\end{center}
		\caption{\change{
		Computation of binary kicks via task-based recursive parallelism. \textbf{Left}: the task of all interactions between two distinct particle sets is divided, by halving each set, and executed in two stages (as indicated) of two mutually independent sub-tasks, which can be done in parallel. \textbf{Right}: the task of all interactions between a set of particles is divided by halving the set. The sub-tasks of interactions within each half are done in parallel first, before the interactions between the two halves are done in a second stage.
		\label{fig:parallelise}
		}}
	\end{minipage}
\end{figure}

\subsection{Efficient c\change{alcul}ation \change{and parallelisation} of binary kicks}
\label{app:eff:Kepler}
\change{The fact that} the composite map $\psi^{W}$ \change{requires exactly the reverse} order \change{of binary} kicks \change{as its adjoint $\psi^{\dag W}$ renders their efficient implementation non-trivial}. \change{Fortunately, these maps are unaffected by a} re-ordering \change{which preserves for each particle the order of its binary-kick interactions. In particular, two} individual maps $\psi^{W_{\ij}}$ and $\psi^{W_{kl}}$ \change{are \emph{mutually independent} and} can be swapped \change{or even computed simultaneously} if all four indices differ. \change{For sufficiently large $N$, this freedom allows} synchronous execution of binary kicks\change{, which} can be implemented by computational parallelism on all levels, including vectorisation.

For \change{vectorisation}, \change{we use} a Kepler solver \change{without} branches \change{(except one to ensure $0\le h<P$ for elliptic orbits)\footnote{\change{The Kepler solver is a modification of one published online as part of another project, whose author obtained it elsewhere on the internet but lost track of its origin. It is based on solving Kepler's equation in universal variables in way that is independent of the nature of the orbit (elliptic, parabolic, or hyperbolic).}}}. Most contemporary \change{CPUs} support vector\change{s of} size \change{$\sub{n}{vec}=4$} for double-precision arithmetic, implying that \change{$\sub{n}{vec}$} Kepler problems can be solved synchronously. \change{A}n efficient way to \change{vectorise} the map $\psi^{W}$ for $N$ particles with $K=N(N-1)/2$ interactions (for the algorithms~\ref{eq:map:HB15} and \change{\ref{eq:map:[DB]^2_4}}) is similar to a round robin sports tournament, where each team plays each other team exactly once. This requires \change{$K/\lfloor N/2\rfloor$} rounds with $\change{\lfloor}N/2\change{\rfloor}$ interactions. \change{As long as $\sub{n}{vec}\le \lfloor N/2\rfloor$, all $K$ interactions can be computed with $\lceil K/\sub{n}{vec}\rceil$ calls to the vectorised Kepler solver, see also Fig.~\ref{fig:vectorise}.

Multi-threaded hardware can be exploited by task-based recursive parallelism using the divide-and-conquer paradigm as explained in Fig.~\ref{fig:parallelise}. To ensure that the order of binary kicks is unaffected by whether or not an interaction task is executed serially or in parallel, the recursive task-based algorithm must also be used with the serial execution down to tasks too small to be split. For the adjoint map, the orders of parallelisation stages and vectorised loops (see Figs.~\ref{fig:vectorise} and \ref{fig:parallelise}) are simply reversed. Remarkably, the requirement of a deterministic order of interactions for each particle renders our parallel implementation deterministic like serial computer programs.

The required ordering of binary-kick interactions implies a rather complex memory access pattern, which in turn hampers computational efficiency as the run-time environment must maintain cache coherence. Despite this, our implementation, \textsc{triton}, achieves good performance and reasonable scalability (strong scaling of 0.72 for 16 cores and $N=1024$). For the [DB]$^2$ simulation reported in Fig.~\ref{fig:nbody}, \textsc{triton} required 13.2 hours on 16 cores or 41 hours on 4 cores, about half of the 3 days reported by \cite{GoncalvesFerrariEtAl2014} for their code \textsc{sakura} (also running on 4 cores with very similar CPU), which also uses a Kepler solver for each particle pair, but only once per time step. Thus, per call to the Kepler solver, \textsc{triton} is almost four times faster than \textsc{sakura}, reflecting the fact that \textsc{triton} is vectorised while \textsc{sakura} is not.
}
\label{lastpage}
\end{document}